\documentclass{amsart}

\usepackage[utf8]{inputenc}
\usepackage[english]{babel}

\usepackage{amscd}
\usepackage{amssymb}
\usepackage{commath}
\usepackage{mathrsfs}
\usepackage{mathtools}

\usepackage{enumitem}
\usepackage{float}
\usepackage{graphicx}
\usepackage{lineno}
\usepackage{multicol}
\usepackage{needspace}
\usepackage{tikz}
\usepackage{tikz-cd}
\usetikzlibrary{decorations.markings}
\usepackage[all]{xy}

\usepackage{hyperref}

\newtheorem{theorem}{Theorem}[section]
\newtheorem{corollary}[theorem]{Corollary}
\newtheorem{lemma}[theorem]{Lemma}
\newtheorem{proposition}[theorem]{Proposition}

\theoremstyle{definition}

\newtheorem{definition}[theorem]{Definition}

\newtheorem{question}[theorem]{Question}
\newtheorem{remark}[theorem]{Remark}

\numberwithin{equation}{section}

\let\mathcal\mathscr

\newcommand{\restr}[2]{#1\!\!\upharpoonright_{#2}}

\begin{document}

%\linenumbers

\title{Colorings and weak separation in graph topologies}

\author{Alejandro R\'ios Herrej\'on}
\address {A. Ríos Herrej\'{o}n\\
Departamento de Matem\'aticas, Facultad de Ciencias, Universidad Nacional Aut\'onoma de M\'exico, Circuito ext. s/n, Ciudad Universitaria, C.P. 04510, Ciudad de M\'exico, M\'exico}
\email{chanchito@ciencias.unam.mx}



\begin{abstract}
In this paper, we study the topological spaces $X_c$ associated with colorings $c : [\omega_1]^2 \to 2$, focusing on the weak separation number and its behavior in products and Pixley--Roy hyperspaces. We characterize the condition $R(X_c) = \omega$ in terms of combinatorial properties of the coloring and construct, for each $n \in \mathbb{N}$, colorings that distinguish the finite levels of these properties. We also introduce several coloring principles and investigate their relationships with other combinatorial principles from the literature. Furthermore, we establish that, for every coloring $c : [\omega_1]^2 \to 2$, we have $c(\mathcal{F}[X_c]) = c^*(\mathcal{F}[X_c]) = R(X_c^\omega)$. Finally, we prove that it is relatively consistent with $\mathsf{ZFC}$ that there exists a Hausdorff space $X_c$ such that $R(X_c) = \omega < c(\mathcal{F}[X_c])$.
\end{abstract}

\maketitle
\enlargethispage{2pt}

%%%%%%%%%%%%%%%%%%%%%%%%%%%%%%%%%%%%%%%%%%%%%%%%%%%%%%%%%%%%%%
\section{Introduction}\label{sec_intro}
%%%%%%%%%%%%%%%%%%%%%%%%%%%%%%%%%%%%%%%%%%%%%%%%%%%%%%%%%%%%%%

The topologies arising from colorings of pairs of ordinals, studied by Hajnal and Juhász \cite{hajju1979}, connect combinatorial features of colorings with cardinal properties of the resulting spaces. We revisit these constructions to investigate the weak separation number, with particular attention to the relationship between the spaces associated with a coloring and its complement.

Combining combinatorial and topological methods, we study separation properties, products, and the role of set-theoretic principles. We also address problems posed in \cite{rios2023} concerning the weak separation number and the cellularity of Pixley--Roy hyperspaces.

The paper is organized as follows. Section~\ref{sec_pre} recalls the necessary preliminaries. Section~\ref{sec_GT} develops the combinatorial framework for graph topologies, while Section~\ref{sec_SGT} investigates their separation properties. Section~\ref{sec_CP} introduces coloring principles and examines their connections with other set-theoretic principles. Section~\ref{sec_PR} studies Pixley--Roy hyperspaces and their associated cardinal invariants. Finally, Section~\ref{sec_Q} presents several open questions.

%%%%%%%%%%%%%%%%%%%%%%%%%%%%%%%%%%%%%%%%%%%%%%%%%%%%%%%%%%%%%%
\section{Preliminaries and background}\label{sec_pre}
%%%%%%%%%%%%%%%%%%%%%%%%%%%%%%%%%%%%%%%%%%%%%%%%%%%%%%%%%%%%%%

In this section, we establish the terminology and notation used throughout the paper and recall some relevant background. We mostly follow the standard conventions in the literature.

The symbol $\omega$ denotes the first infinite cardinal, and $\mathbb{N} := \omega \setminus \{0\}$. We also write $\mathfrak{c}$ for the cardinal $2^\omega$. For any set $X$ and any positive cardinal $\kappa$, we use the notation $[X]^{\leq\kappa}$, $[X]^{<\kappa}$, and $[X]^\kappa$ for the collections of subsets of $X$ of cardinality at most $\kappa$, less than $\kappa$, and exactly $\kappa$, respectively.

Let $X$ be a topological space, and let $\tau_X$ denote its topology. For any topological cardinal function $\varphi$, we define
\[
\varphi^*(X) := \sup\{\varphi(X^n) : n \in \mathbb{N}\}.
\]
If $I$ is a nonempty set, an \emph{$I$-assignment} is a set $\{(x_\alpha,U_\alpha) : \alpha \in I\} \subseteq X \times \tau_X$ such that $x_\alpha \in U_\alpha$ for every $\alpha \in I$. When the index set is clear from the context, we simply use the term \emph{assignment}.

A space $X$ is \emph{weakly separated} if there exists an assignment $\{(x,U_x) : x \in X\}$ such that, for every pair of distinct points $x,y \in X$, either $x \not\in U_y$ or $y \not\in U_x$. A family $\{U_x : x \in X\}$ satisfying these conditions is called a \emph{weak separation} of $X$. Weakly separated spaces were introduced by Tkachenko in \cite{tkachenko1978}. The \emph{weak separation number} of $X$ is the cardinal
\[
R(X) := \sup\{|A| : A \text{ is a weakly separated subspace of } X\} + \omega.
\]
The condition $R(X) = \omega$ corresponds to the property $\mathrm{HG}$ of Hart and Kunen (see \cite{harkun2020}), while $R(X^\omega) = \omega$ corresponds to its stronger version $\mathrm{stHG}$.

For any topological space $X$, let $\mathcal{F}(X) := [X]^{<\omega} \setminus \{\emptyset\}$. Moreover, for any $F \in \mathcal{F}(X)$ and $S \subseteq X$, define
\[
[F,S] := \{G \in \mathcal{F}(X) : F \subseteq G \subseteq S\}.
\]
The collection $\{[F,U] : F \in \mathcal{F}(X) \ \wedge\ U \in \tau_X\}$ forms a base for a topology on $\mathcal{F}(X)$, known as the \emph{Pixley--Roy topology}. We denote by $\mathcal{F}[X]$ the set $\mathcal{F}(X)$ equipped with this topology.

We now recall some results concerning the weak separation number, products, and Pixley--Roy hyperspaces. It was shown in \cite{rios2023} that, for every topological space $X$, the following relations hold:
\begin{equation}\label{eq_1}
R(X) \leq c(\mathcal{F}[X]) \leq c^*(\mathcal{F}[X]) = R(X^\omega) = R^*(X) \leq nw(X).
\end{equation}
In particular, since $\mathsf{MA}(\omega_1)$ implies that every space of countable cellularity has precaliber $\omega_1$, and this property is preserved under finite products, the preceding relations imply that, under $\mathsf{MA}(\omega_1)$, if $c(\mathcal{F}[X]) = \omega$, then $R(X^\omega) = \omega$.

On the other hand, Juhász, Soukup, and Szentmiklóssy proved that, for each $m \in \mathbb{N}$, it is consistent that there exists a space $X$ such that $R(X^m) = \omega < R(X^{m+1})$ (see \cite[Theorem~3.1]{jss1996}). Furthermore, Gruenhage explained how to obtain, from a $2$-entangled set of real numbers, a space $X$ such that $R(X) = \omega$ and $X^2$ contains an uncountable discrete subspace (see \cite[p.~304]{gruenhage1989}). Finally, since the existence of $2$-entangled sets is consistent with $\mathsf{MA}(\omega_1)$ (see \cite[Theorem~97]{guztod2024}), this axiom does not guarantee that the property $R(X) = \omega$ is preserved under finite products.

%%%%%%%%%%%%%%%%%%%%%%%%%%%%%%%%%%%%%%%%%%%%%%%%%%%%%%%%%%%%%%
\section{Graph topology}\label{sec_GT}
%%%%%%%%%%%%%%%%%%%%%%%%%%%%%%%%%%%%%%%%%%%%%%%%%%%%%%%%%%%%%%

Let $\mathrm{Fn}(\omega_1,2)$ be the set of all finite partial functions from $\omega_1$ to $2$, that is, $p \in \mathrm{Fn}(\omega_1,2)$ if and only if $p \in [\omega_1 \times 2]^{<\omega}$ and $p$ is a function.

Let $c : [\omega_1]^2 \to 2$ be a coloring. From now on, we adopt the convention that $c(\{\alpha\}) := 0$ for every $\alpha < \omega_1$. The \emph{complementary coloring} of $c$ is the coloring $\overline{c} : [\omega_1]^2 \to 2$ defined by
\[
\overline{c}(\{\alpha,\beta\}) := 1 - c(\{\alpha,\beta\}).
\]
Note that $\overline{c}$ also satisfies the preceding convention, that is, $\overline{c}(\{\alpha\}) := 0$ for every $\alpha < \omega_1$.

For each $i < 2$, the symbol $c^{-1}\{i\}$ denotes the \emph{fiber} of $i$ under $c$, that is,
\[
c^{-1}\{i\} := \big\{\{\alpha,\beta\} \in [\omega_1]^2 : c(\{\alpha,\beta\}) = i\big\}.
\]
Moreover, for every $\alpha < \omega_1$ and $i < 2$, let
\[
U_\alpha^i(c) := \{\beta < \omega_1 : c(\{\alpha,\beta\}) = i\}.
\]
Let $\tau_c$ be the topology on $\omega_1$ generated by the subbase $\{U_\alpha^i(c) : \alpha < \omega_1 \ \wedge \ i < 2\}$. Finally, let $U_\emptyset(c) := \omega_1$, and for each $p \in \mathrm{Fn}(\omega_1,2) \setminus \{\emptyset\}$, let
\[
\textstyle U_p(c) := \bigcap\big\{U_\alpha^{p(\alpha)}(c) : \alpha \in \operatorname{dom}(p)\big\}.
\]
When the coloring $c$ is clear from the context, we omit it from the notation and simply write $U_\alpha^i$ and $U_p$ instead of $U_\alpha^i(c)$ and $U_p(c)$, respectively.

In \cite{hajju1979}, Hajnal and Juhász explained that the resulting topological space, $X_c := (\omega_1,\tau_c)$, is zero-dimensional and that the collection $\{U_p : p \in \mathrm{Fn}(\omega_1,2)\}$ forms a base of clopen sets.

In Definition~\ref{def_GT_central}, we introduce some notions that will be fundamental to the development of this work.

\begin{definition}\label{def_GT_central}
Let $n \in \mathbb{N}$ and let $c : [\omega_1]^2 \to 2$ be a coloring.

\begin{enumerate}[label=(\arabic*)]

\item For each $\alpha < \omega_1$, let $c_\alpha : \omega_1 \to 2$ be the function defined by
\[
c_\alpha(\gamma) := c(\{\alpha,\gamma\}).
\]

\item For every distinct $\alpha,\beta < \omega_1$, define
\[
\Delta_c(\alpha,\beta) :=
\{\gamma \in \omega_1 \setminus \{\alpha,\beta\} : c_\alpha(\gamma) \neq c_\beta(\gamma)\}.
\]

\item For each $\alpha < \omega_1$, let
\[
\mathcal{P}_\alpha := \big\{p \in \mathrm{Fn}(\omega_1,2) \setminus \{\emptyset\} : p = \restr{c_\alpha}{\operatorname{dom}(p)}\big\}
\]
and $\mathcal{B}_\alpha := \{U_p : p \in \mathcal{P}_\alpha\}$.

\item For any $J \in [\omega_1]^{\omega_1}$, we say that a family $\{E_\alpha : \alpha \in J\} \subseteq [\omega_1]^{<\omega} \setminus \{\emptyset\}$ is a \emph{block sequence} if $\max(E_\alpha) < \min(E_\beta)$ whenever $\alpha,\beta \in J$ and $\alpha < \beta$. If, in addition, $|E_\alpha| = n$ for every $\alpha \in J$, we say that the family is an \emph{$n$-block sequence}.

\item We say that condition $(\star_n^c)$ holds if, for every $J \in [\omega_1]^{\omega_1}$ and every $n$-block sequence $\{E_\alpha : \alpha \in J\}$, there exist distinct $\alpha,\beta \in J$ such that
\[
\restr{c_\alpha}{E_\alpha \cup E_\beta} = \restr{c_\beta}{E_\alpha \cup E_\beta}.
\]

When the coloring $c$ is clear from the context, we simply write $(\star_n)$ instead of $(\star_n^c)$. We also define
\[
S_c := \{m \in \mathbb{N} : \text{condition $(\star_m)$ holds}\}.
\]

\end{enumerate}

\end{definition}

Observe that if $c : [\omega_1]^2 \to 2$ is a coloring and $\alpha,\beta < \omega_1$ are distinct, then
\[
\Delta_c(\alpha,\beta) = \Delta_c(\beta,\alpha) = \Delta_{\overline{c}}(\alpha,\beta).
\]
Moreover, for each $\alpha < \omega_1$, the collection $\mathcal{B}_\alpha$ is a local base at $\alpha$ in $X_c$.

\begin{remark}\label{rmk_GT_base_local} If $\alpha,\beta < \omega_1$ and $q \in \mathcal{P}_\beta$, then $\alpha \in U_q$ if and only if $\restr{c_\alpha}{\operatorname{dom}(q)} = \restr{c_\beta}{\operatorname{dom}(q)}$. Consequently, for every $\alpha,\beta < \omega_1$, $p \in \mathcal{P}_\alpha$, and $q \in \mathcal{P}_\beta$, the following statements are equivalent:

\begin{enumerate}
\item $\alpha \in U_q$ and $\beta \in U_p$.

\item $\restr{c_\alpha}{\operatorname{dom}(p) \cup \operatorname{dom}(q)} = \restr{c_\beta}{\operatorname{dom}(p) \cup \operatorname{dom}(q)}$.
\end{enumerate}

\end{remark}

The first result of this section highlights the importance of condition $(\star_n)$ for our purposes:

\begin{lemma}\label{lema_GT_Xc_HG} The following statements are equivalent for any coloring $c : [\omega_1]^2 \to 2$.

\begin{enumerate}
\item $R(X_c) = \omega$.

\item For every $J \in [\omega_1]^{\omega_1}$ and every family $\{F_\alpha : \alpha \in J\} \subseteq [\omega_1]^{<\omega} \setminus \{\emptyset\}$, there exist distinct $\alpha,\beta \in J$ such that $\restr{c_\alpha}{F_\alpha \cup F_\beta} = \restr{c_\beta}{F_\alpha \cup F_\beta}$.

\item $S_c = \mathbb{N}$.

\end{enumerate}

\end{lemma}

\begin{proof} Our goal is to show that (1) and (2) are equivalent, as are (2) and (3).

To prove that (1) implies (2), suppose that $R(X_c) = \omega$. Let $J \in [\omega_1]^{\omega_1}$ and let $\{F_\alpha : \alpha \in J\} \subseteq [\omega_1]^{<\omega} \setminus \{\emptyset\}$. Moreover, for each $\alpha \in J$, let $p_\alpha := \restr{c_\alpha}{F_\alpha}$. Then $\{(\alpha,U_{p_\alpha}) : \alpha \in J\}$ is an assignment. Hence, there exist distinct $\alpha,\beta \in J$ such that $\alpha \in U_{p_\beta}$ and $\beta \in U_{p_\alpha}$. Consequently, by Remark~\ref{rmk_GT_base_local}, we have $\restr{c_\alpha}{F_\alpha \cup F_\beta} = \restr{c_\beta}{F_\alpha \cup F_\beta}$.

To prove that (2) implies (1), suppose that $R(X_c) > \omega$. Let $J \in [\omega_1]^{\omega_1}$ be a weakly separated subspace of $X_c$ and, for each $\alpha \in J$, let $p_\alpha \in \mathcal{P}_\alpha$ be such that $\{U_{p_\alpha} : \alpha \in J\}$ is a weak separation of $J$. If we define $F_\alpha := \operatorname{dom}(p_\alpha)$ for each $\alpha \in J$, then $\{F_\alpha : \alpha \in J\} \subseteq [\omega_1]^{<\omega} \setminus \{\emptyset\}$. Moreover, if $\alpha,\beta \in J$ are distinct and we assume, without loss of generality, that $\alpha \notin U_{p_\beta}$, then $\restr{c_\alpha}{F_\beta} \neq \restr{c_\beta}{F_\beta}$. Thus, $\restr{c_\alpha}{F_\alpha \cup F_\beta} \neq \restr{c_\beta}{F_\alpha \cup F_\beta}$, contradicting (2).

The implication from (2) to (3) is immediate. To prove that (3) implies (2), let $J \in [\omega_1]^{\omega_1}$ and let $\{F_\alpha : \alpha \in J\} \subseteq [\omega_1]^{<\omega} \setminus \{\emptyset\}$. For each $\alpha \in J$, define $H_\alpha := F_\alpha \cup \{\alpha\}$. By the $\Delta$-system lemma, there exist $I \in [J]^{\omega_1}$ and $H \in [\omega_1]^{<\omega}$ such that $H_\alpha \cap H_\beta = H$ for every pair of distinct $\alpha,\beta \in I$.

Moreover, there exists $K \in [I]^{\omega_1}$ such that $K \cap H = \emptyset$ and $\restr{c_\alpha}{H} = \restr{c_\beta}{H}$ for every $\alpha,\beta \in K$. Now, for each $\alpha \in K$, define $E_\alpha := H_\alpha \setminus H$. Since the function $K \to \mathbb{N}$ given by $\alpha \mapsto |E_\alpha|$ has a fiber of cardinality $\omega_1$, there exist $L \in [K]^{\omega_1}$ and $n \in \mathbb{N}$ such that $|E_\alpha| = n$ for every $\alpha \in L$.

\medskip

\noindent {\bf Claim.} There exists $M \in [L]^{\omega_1}$ such that $\max(E_\alpha) < \min(E_\beta)$ whenever $\alpha,\beta \in M$ satisfy $\alpha < \beta$.

\medskip

The strategy is to construct, by transfinite recursion on $\omega_1$, a sequence $\{\alpha_\xi : \xi < \omega_1\} \subseteq L$ such that, whenever $\xi < \eta < \omega_1$, we have $\alpha_\xi < \alpha_\eta$ and $\max(E_{\alpha_\xi}) < \min(E_{\alpha_\eta})$. Suppose that, for some $\xi < \omega_1$, we have constructed the family $\{\alpha_\eta : \eta < \xi\}$ with the desired properties. Since $\delta := \sup\{\max(E_{\alpha_\eta}) : \eta < \xi\}+1$ is a countable ordinal, the set $\{\alpha \in L : E_\alpha \cap \delta \neq \emptyset\}$ is countable, because $\{E_\alpha : \alpha \in L\}$ is a pairwise disjoint family. Moreover, $\sup\{\alpha_\eta : \eta < \xi\}<\omega_1$. Therefore, there exists $\alpha_\xi \in L$ such that $\alpha_\xi > \sup\{\alpha_\eta : \eta < \xi\}$ and $E_{\alpha_\xi} \cap \delta = \emptyset$. Consequently, $\alpha_\eta < \alpha_\xi$ and $\max(E_{\alpha_\eta}) < \min(E_{\alpha_\xi})$ for every $\eta < \xi$. This completes the recursion.

\medskip

Now, since $\{E_\alpha : \alpha \in M\}$ is an $n$-block sequence, the fact that $n \in S_c$ guarantees the existence of distinct $\alpha,\beta \in M$ such that $\restr{c_\alpha}{E_\alpha \cup E_\beta} = \restr{c_\beta}{E_\alpha \cup E_\beta}$. Finally, since we also have $\restr{c_\alpha}{H} = \restr{c_\beta}{H}$ and $F_\alpha \cup F_\beta \subseteq H \cup E_\alpha \cup E_\beta$, we conclude that $\restr{c_\alpha}{F_\alpha \cup F_\beta} = \restr{c_\beta}{F_\alpha \cup F_\beta}$.
\end{proof}

Next, we study some basic properties of condition $(\star_n)$. We begin with the following result:

\begin{lemma}\label{lema_GT_estrella_baja} If $m,n \in \mathbb{N}$ satisfy $m > n$, $c : [\omega_1]^2 \to 2$ is a coloring, and $(\star_m)$ holds, then $(\star_n)$ also holds.

\end{lemma}

\begin{proof} Let $J \in [\omega_1]^{\omega_1}$ and let $\{E_\alpha : \alpha \in J\}$ be an $n$-block sequence. We can find $I \in [J]^{\omega_1}$ such that $\max(E_\alpha)+(m-n) < \min(E_\beta)$ whenever $\alpha,\beta \in I$ and $\alpha < \beta$. Now, if for each $\alpha \in I$ we define
\[
F_\alpha := E_\alpha \cup \{\max(E_\alpha)+i+1 : i < m-n\},
\]
then $\{F_\alpha : \alpha \in I\}$ is an $m$-block sequence. Therefore, condition $(\star_m)$ guarantees the existence of distinct $\alpha,\beta \in I$ such that $\restr{c_\alpha}{F_\alpha \cup F_\beta} = \restr{c_\beta}{F_\alpha \cup F_\beta}$. In particular, since $E_\alpha \cup E_\beta \subseteq F_\alpha \cup F_\beta$, we have $\restr{c_\alpha}{E_\alpha \cup E_\beta} = \restr{c_\beta}{E_\alpha \cup E_\beta}$, which proves that condition $(\star_n)$ holds.
\end{proof}

Lemma~\ref{lema_GT_estrella_baja} implies that, for any coloring $c : [\omega_1]^2 \to 2$, the set $S_c$ is an initial segment of $\mathbb{N}$. Therefore, there are only three possibilities: $S_c = \emptyset$, $S_c = \{1,\ldots,n\}$ for some $n \in \mathbb{N}$, or $S_c = \mathbb{N}$.

Now, if $\underline{0} : [\omega_1]^2 \to 2$ is the constant zero coloring and $\underline{1} : [\omega_1]^2 \to 2$ is the constant one coloring, it is immediate that $S_{\underline{0}} = \mathbb{N}$ and $S_{\underline{1}} = \emptyset$. Thus, to show that all three possibilities for $S_c$ can occur, it remains only to find colorings $c : [\omega_1]^2 \to 2$ for which $S_c$ is a finite nonempty initial segment of $\mathbb{N}$.

To achieve this goal, we first establish some preliminary results that will simplify the subsequent proofs. In this regard, Remark~\ref{rmk_GT_particion} and Lemmas~\ref{lema_GT_adelgazamiento_1} and \ref{lema_GT_adelgazamiento_2} play a fundamental role.

\begin{remark}\label{rmk_GT_particion} An important fact, which can be proved by transfinite induction, is that if $n \in \mathbb{N}$ and $\{B_\xi : \xi < \omega_1\}$ is an $(n+2)$-block sequence that partitions $\omega_1$, then
\[
B_\xi = \{(n+2)\xi + i : i < n+2\}
\]
for every $\xi < \omega_1$. Therefore, for the unique family $\{B_\xi : \xi < \omega_1\}$ with these properties, we can define, for every $\xi < \omega_1$ and $i < n+1$,
\[
a_\xi := (n+2)\xi \quad \text{and} \quad b_\xi^i := a_\xi + i + 1.
\]
Thus, $B_\xi = \{a_\xi\} \cup \{b_\xi^i : i < n+1\}$ for every $\xi < \omega_1$. We will use this notation in the proofs of Lemma~\ref{lema_GT_adelgazamiento_1} and Proposition~\ref{prop_GT_estrella_no_sube}.

\end{remark}

\begin{lemma}\label{lema_GT_adelgazamiento_1} Let $n \in \mathbb{N}$, let $c : [\omega_1]^2 \to n+1$ be a coloring, and let $\{B_\xi : \xi < \omega_1\}$ be the $(n+2)$-block sequence that partitions $\omega_1$. Moreover, let $J \in [\omega_1]^{\omega_1}$ and let $\{E_\alpha : \alpha \in J\}$ be an $n$-block sequence. For each $\gamma < \omega_1$, let $\rho(\gamma)$ be the unique ordinal $\xi < \omega_1$ such that $\gamma \in B_\xi$. Also, for each $\alpha \in J$, let
\[
T_\alpha := \{\rho(\alpha)\} \cup \{\rho(\gamma) : \gamma \in E_\alpha\}.
\]

Under these assumptions, there exist $m \in \mathbb{N}$ with $m \leq n+1$, $r < m$, a family of colors $\{k_s : s \in m \setminus \{r\}\} \subseteq n+1$, a set $D \in [n+1]^{\leq n}$, and a set $I \in [J]^{\omega_1}$ such that the following properties hold:

\begin{enumerate}[label=(\arabic*)]
\item The family $\{T_\alpha : \alpha \in I\}$ is pairwise disjoint, and $|T_\alpha| = m$ for every $\alpha \in I$.

\item If, for each $\alpha \in I$, $\{t_\alpha^s : s < m\}$ is the strictly increasing enumeration of $T_\alpha$, then $\rho(\alpha) = t_\alpha^r$ for every $\alpha \in I$.

\item $c(\{t_\alpha^r,t_\alpha^s\}) = k_s$ for every $\alpha \in I$ and $s \in m \setminus \{r\}$.

\item Either $\alpha = a_{t_\alpha^r}$ for every $\alpha \in I$, or there exists $j < n+1$ such that $\alpha = b_{t_\alpha^r}^j$ for every $\alpha \in I$.

\item In the first case of (4), we have
\[
D = \{i < n+1 : b_{t_\alpha^r}^i \in E_\alpha\}
\]
for every $\alpha \in I$. In the second case, $D = \{j\}$.

\end{enumerate}

\end{lemma}

\begin{proof} First, for each $\xi < \omega_1$, the set $J_\xi := \{\alpha \in J : \xi \in T_\alpha\}$ is finite. Indeed, if $\alpha \in J_\xi$, then $\alpha \in B_\xi$ or $E_\alpha \cap B_\xi \neq \emptyset$. The first possibility occurs for at most $n+2$ indices, and the same holds for the second, since the sets $E_\alpha$ are pairwise disjoint and $|B_\xi| = n+2$. Therefore, $|J_\xi| \leq 2(n+2)$ for every $\xi < \omega_1$. Consequently, for each $\alpha \in J$, the set
\[
\textstyle
\{\beta \in J : T_\alpha \cap T_\beta \neq \emptyset\}
= \bigcup_{\xi \in T_\alpha} J_\xi
\]
is finite.

This allows us to construct, by transfinite recursion, a sequence $\{\alpha_\eta : \eta < \omega_1\}$ of distinct elements of $J$ such that the family $\{T_{\alpha_\eta} : \eta < \omega_1\}$ is pairwise disjoint. Indeed, at each stage $\eta < \omega_1$, it suffices to avoid a countable union of finite sets, while $J$ is uncountable.

Thus, the set $I_0 := \{\alpha_\eta : \eta < \omega_1\}$ satisfies $I_0 \in [J]^{\omega_1}$, and the family $\{T_\alpha : \alpha \in I_0\}$ is pairwise disjoint. Moreover, since $|E_\alpha| = n$, we have $1 \leq |T_\alpha| \leq n+1$ for every $\alpha \in I_0$. Consequently, there exist $m \in \mathbb{N}$, with $m \leq n+1$, and $I_1 \in [I_0]^{\omega_1}$ such that $|T_\alpha| = m$ for every $\alpha \in I_1$.

For each $\alpha \in I_1$, let $\{t_\alpha^s : s < m\}$ be the strictly increasing enumeration of $T_\alpha$. Since $\rho(\alpha) \in T_\alpha$, there exists a unique $r_\alpha < m$ such that $\rho(\alpha) = t_\alpha^{r_\alpha}$. As there are only $m$ possibilities for $r_\alpha$, there exist $r < m$ and $I_2 \in [I_1]^{\omega_1}$ such that $\rho(\alpha) = t_\alpha^r$ for every $\alpha \in I_2$.

On the other hand, for each $\alpha \in I_2$, consider the function $f_\alpha : m \setminus \{r\} \to n+1$ defined by
\[
f_\alpha(s) := c(\{t_\alpha^r, t_\alpha^s\}).
\]
Since both the domain and the codomain are finite, there are only finitely many possibilities for $f_\alpha$. Therefore, there exist a function $f : m \setminus \{r\} \to n+1$ and a set $I_3 \in [I_2]^{\omega_1}$ such that $f_\alpha = f$ for every $\alpha \in I_3$. Thus, setting $k_s := f(s)$ for each $s \in m \setminus \{r\}$, we obtain $c(\{t_\alpha^r, t_\alpha^s\}) = k_s$ for every $\alpha \in I_3$ and $s \in m \setminus \{r\}$.

Now, since $\rho(\alpha) = t_\alpha^r$ for every $\alpha \in I_3$, it follows that
\[
\alpha \in B_{t_\alpha^r}
= \{a_{t_\alpha^r}\} \cup \{b_{t_\alpha^r}^j : j < n+1\}.
\]
Since there are only $n+2$ possibilities for the position of $\alpha$ within its block, we can find $I_4 \in [I_3]^{\omega_1}$ such that one of the following alternatives holds:

\begin{itemize}
\item $\alpha = a_{t_\alpha^r}$ for every $\alpha \in I_4$.
\item There exists $j < n+1$ such that $\alpha = b_{t_\alpha^r}^j$ for every $\alpha \in I_4$.
\end{itemize}

If the second alternative holds, it suffices to take $I := I_4$ and $D := \{j\}$. If the first alternative holds, for each $\alpha \in I_4$, define
\[
D_\alpha := \{i < n+1 : b_{t_\alpha^r}^i \in E_\alpha\}.
\]
Observe that, for each $\alpha \in I_4$, the condition $|E_\alpha| = n$ implies that $|D_\alpha| \leq n$. Moreover, since $D_\alpha \subseteq n+1$, there are only finitely many possibilities for $D_\alpha$. Consequently, there exist $D \in [n+1]^{\leq n}$ and $I \in [I_4]^{\omega_1}$ such that $D_\alpha = D$ for every $\alpha \in I$.

In either case, we have obtained a set $I \in [J]^{\omega_1}$, the parameters $m,r$, the colors $k_s$, and the set $D$ satisfying all five conditions of the statement. This completes the proof.
\end{proof}

The proof of Lemma~\ref{lema_GT_adelgazamiento_2} follows the same technique as that of Lemma~\ref{lema_GT_adelgazamiento_1}. Since the argument is simpler and requires no additional ideas, we omit the details.

\begin{lemma}\label{lema_GT_adelgazamiento_2} Let $n \in \mathbb{N}$, let $c : [\omega_1]^2 \to 2$ be a coloring, let $J \in [\omega_1]^{\omega_1}$, and let $\{E_\alpha : \alpha \in J\}$ be an $n$-block sequence. If, for each $\alpha \in J$, we set $T_\alpha := E_\alpha \cup \{\alpha\}$, then there exist $m \in \mathbb{N}$ with $m \leq n+1$, $r < m$, a family of colors $\{k_s : s \in m \setminus \{r\}\} \subseteq 2$, and a set $I \in [J]^{\omega_1}$ such that the following properties hold:

\begin{enumerate}[label=(\arabic*)]
\item The family $\{T_\alpha : \alpha \in I\}$ is pairwise disjoint, and $|T_\alpha| = m$ for every $\alpha \in I$.

\item If, for each $\alpha \in I$, $\{t_\alpha^s : s < m\}$ is the strictly increasing enumeration of $T_\alpha$, then $\alpha = t_\alpha^r$ for every $\alpha \in I$.

\item $c(\{t_\alpha^r,t_\alpha^s\}) = k_s$ for every $\alpha \in I$ and $s \in m \setminus \{r\}$.
\end{enumerate}

\end{lemma}

The final goal of this section is to construct, in Proposition~\ref{prop_GT_estrella_no_sube}, for each $n \in \mathbb{N}$, a coloring $c_n : [\omega_1]^2 \to 2$ such that $S_{c_n} = \{1,\ldots,n\}$ or, equivalently, such that $(\star_n)$ holds but $(\star_{n+1})$ does not. We adopt the convention that, if $m < \omega$ and $a \in [\omega_1]^m$, then $a(i)$ denotes the $i$-th element of $a$ in its strictly increasing enumeration. To achieve our goal, we use the following result of Peng and Wu~\cite[Corollary~3, p.~218]{penwu2018}:

\begin{theorem}\label{thm_GT_PW} The principle $\operatorname{Pr}_0(\omega_1,\omega_1,n)$ holds for every $n<\omega$. That is, for each $n<\omega$, there exists a coloring $c:[\omega_1]^2\to\omega_1$ such that, for every $m<n$, every pairwise disjoint family $\{a_\alpha:\alpha<\omega_1\}$ of elements of $[\omega_1]^m$, and every function $M:m\times m\to\omega_1$, there exist $\alpha<\beta<\omega_1$ such that
\[
c(\{a_\alpha(i),a_\beta(j)\})=M(i,j)
\]
for every $i,j<m$.

\end{theorem}

Our first task is to reduce the number of colors in the coloring provided by Theorem~\ref{thm_GT_PW} to a finite number while preserving its combinatorial properties. To this end, Lemma~\ref{lema_GT_PW_n+1} will be essential for proving Proposition~\ref{prop_GT_estrella_no_sube}.

\begin{lemma}\label{lema_GT_PW_n+1} For each $n \in \mathbb{N}$, there exists a coloring $\widehat{c_n} : [\omega_1]^2 \to n+1$ such that, for every $m < n+2$, every pairwise disjoint family $\{a_\alpha : \alpha < \omega_1\}$ of elements of $[\omega_1]^m$, and every function $M : m\times m \to n+1$, there exist $\alpha < \beta < \omega_1$ such that
\[
\widehat{c_n}(\{a_\alpha(i),a_\beta(j)\}) = M(i,j)
\]
for every $i,j < m$.

\end{lemma}

\begin{proof} Let $n \in \mathbb{N}$, and let $d_n : [\omega_1]^2 \to \omega_1$ be the coloring provided by Theorem~\ref{thm_GT_PW} upon applying the principle $\operatorname{Pr}_0(\omega_1,\omega_1,n+2)$. Moreover, let $h_n : \omega_1 \to n+1$ be the function defined by
\[
h_n(\xi) :=
\begin{cases}
\xi, & \text{if } \xi < n+1,\\
0, & \text{if } \xi \geq n+1.
\end{cases}
\]
Clearly, the coloring $\widehat{c_n} := h_n \circ d_n$ satisfies the desired property.
\end{proof}

\begin{proposition}\label{prop_GT_estrella_no_sube} For each $n \in \mathbb{N}$, there exists a coloring $c_n : [\omega_1]^2 \to 2$ such that $(\star_n)$ holds but $(\star_{n+1})$ does not. In particular, $S_{c_n} = \{1,\ldots,n\}$.

\end{proposition}

\begin{proof} Let $n \in \mathbb{N}$, let $\{B_\xi : \xi < \omega_1\}$ be the $(n+2)$-block sequence that partitions $\omega_1$ (see Remark~\ref{rmk_GT_particion}), let $\widehat{c_n} : [\omega_1]^2 \to n+1$ be the coloring provided by Lemma~\ref{lema_GT_PW_n+1}, and let $c_n : [\omega_1]^2 \to 2$ be the coloring defined by
\[
c_n(\{a_\xi,b_\eta^i\}) :=
\begin{cases}
1, & \text{if } \xi \neq \eta \text{ and } \widehat{c_n}(\{\xi,\eta\}) = i,\\
0, & \text{otherwise}.
\end{cases}
\]
We also assign the value $0$ to all pairs not of the form $\{a_\xi,b_\eta^i\}$.

Observe that if $\xi < \omega_1$ and $i < n+1$, then $c_n(\{a_\xi,b_\xi^i\}) = 0$. Moreover, if $\xi < \eta < \omega_1$, there is exactly one $i < n+1$ such that $c_n(\{a_\xi,b_\eta^i\}) = 1$, namely, $i = \widehat{c_n}(\{\xi,\eta\})$. The remainder of the proof is divided into two claims.

\medskip

\noindent {\bf Claim 1.} $c_n$ does not satisfy condition $(\star_{n+1})$.

\medskip

Let $J:=\{a_\xi:\xi<\omega_1\}$ and, for each $\xi<\omega_1$, define $E_{a_\xi}:=\{b_\xi^i:i<n+1\}$. Since $\{B_\xi:\xi<\omega_1\}$ is an $(n+2)$-block sequence, it follows that $\{E_\alpha:\alpha\in J\}$ is an $(n+1)$-block sequence. Now, let $\xi<\eta<\omega_1$ and set $i:=\widehat{c_n}(\{\xi,\eta\})$. By the definition of $c_n$, we have
\[
c_n(\{a_\xi,b_\xi^i\})=0
\quad\text{and}\quad
c_n(\{a_\eta,b_\xi^i\})=1.
\]
Therefore, since $b_\xi^i\in E_{a_\xi}$, we conclude that $c_n$ does not satisfy condition $(\star_{n+1})$.

\medskip

\noindent {\bf Claim 2.} $c_n$ satisfies condition $(\star_n)$.

\medskip

Let $J \in [\omega_1]^{\omega_1}$ and let $\{E_\alpha : \alpha \in J\}$ be an $n$-block sequence. By Lemma~\ref{lema_GT_adelgazamiento_1}, there exist $m \in \mathbb N$, with $m \leq n+1$, $r < m$, a family of colors $\{k_s : s \in m \setminus \{r\}\} \subseteq n+1$, a set $D \in [n+1]^{\leq n}$, and a set $I \in [J]^{\omega_1}$ satisfying the five conditions of that lemma. Moreover, since $|D| \leq n$, we can choose a color $k \in (n+1) \setminus D$.

Let $M : m \times m \to n+1$ be the function defined by
\[
M(s,t) :=
\begin{cases}
k, & \text{if } s = t = r,\\
k_s, & \text{if } t = r \text{ and } s \neq r,\\
k_t, & \text{if } s = r \text{ and } t \neq r,\\
0, & \text{otherwise}.
\end{cases}
\]
Since $\{T_\alpha : \alpha \in I\}$ is a pairwise disjoint family of elements of $[\omega_1]^m$, Lemma~\ref{lema_GT_PW_n+1} guarantees the existence of distinct $\alpha,\beta \in I$ such that
\[
\widehat{c_n}(\{t_\alpha^p, t_\beta^q\}) = M(p,q)
\]
for every $p,q < m$.

Let $\xi := t_\alpha^r$ and $\eta := t_\beta^r$. Since $T_\alpha \cap T_\beta = \emptyset$, we have $\xi \neq \eta$. Moreover, $\widehat{c_n}(\{\xi,\eta\}) = \widehat{c_n}(\{t_\alpha^r, t_\beta^r\}) = M(r,r) = k$. Also, by condition (3) of Lemma~\ref{lema_GT_adelgazamiento_1} and the definition of $M$, for each $s \in m \setminus \{r\}$, we have
\[
\widehat{c_n}(\{\xi, t_\alpha^s\})
= \widehat{c_n}(\{\eta, t_\beta^s\}) = k_s
= M(r,s)
= \widehat{c_n}(\{\xi, t_\beta^s\})
= \widehat{c_n}(\{\eta, t_\alpha^s\}).
\]
Consequently,
\begin{equation}\label{eq_2}
\widehat{c_n}(\{\xi,\zeta\}) = \widehat{c_n}(\{\eta,\zeta\})
\end{equation}
for every $\zeta \in (T_\alpha \cup T_\beta) \setminus \{\xi,\eta\}$.

Our next goal is to show that $c_n(\{\alpha,\gamma\}) = c_n(\{\beta,\gamma\})$ for every $\gamma \in E_\alpha \cup E_\beta$. Let $\gamma \in E_\alpha \cup E_\beta$ and $\zeta := \rho(\gamma)$. Observe that $\zeta \in T_\alpha \cup T_\beta$. To establish the desired equality, we distinguish the two cases described in condition (4) of Lemma~\ref{lema_GT_adelgazamiento_1}.

\medskip

\noindent {\bf Case 1.} $\alpha = a_\xi$ and $\beta = a_\eta$.

\medskip

If $\gamma = a_\zeta$, then $c_n(\{a_\xi, a_\zeta\}) = 0 = c_n(\{a_\eta, a_\zeta\})$. Thus, suppose that $\gamma = b_\zeta^i$ for some $i < n+1$. If $\zeta \notin \{\xi,\eta\}$, equality~\eqref{eq_2} implies that $c_n(\{a_\xi, b_\zeta^i\}) = c_n(\{a_\eta, b_\zeta^i\})$. On the other hand, if $\zeta = \xi$, then $\gamma \in E_\alpha$, since $T_\alpha \cap T_\beta = \emptyset$. Moreover, by condition (5) of Lemma~\ref{lema_GT_adelgazamiento_1}, we have $i \in D$. Since $\widehat{c_n}(\{\eta,\xi\}) = k \notin D$, it follows that $c_n(\{a_\xi, b_\xi^i\}) = 0 = c_n(\{a_\eta, b_\xi^i\})$. Finally, if $\zeta = \eta$, a similar argument shows that $c_n(\{a_\xi, b_\eta^i\}) = 0 = c_n(\{a_\eta, b_\eta^i\})$.

\medskip

\noindent {\bf Case 2.} There exists $j < n+1$ such that $\alpha = b_\xi^j$ and $\beta = b_\eta^j$.

\medskip

In this case, by condition (5) of Lemma~\ref{lema_GT_adelgazamiento_1}, we have $D = \{j\}$; in particular, $\widehat{c_n}(\{\xi,\eta\}) = k \neq j$. Now, if $\gamma = b_\zeta^i$ for some $i < n+1$, then $c_n(\{b_\xi^j, b_\zeta^i\}) = 0 = c_n(\{b_\eta^j, b_\zeta^i\})$. Thus, suppose that $\gamma = a_\zeta$. If $\zeta \notin \{\xi,\eta\}$, equality~\eqref{eq_2} implies that $c_n(\{b_\xi^j, a_\zeta\}) = c_n(\{b_\eta^j, a_\zeta\})$. On the other hand, if $\zeta = \xi$, then $c_n(\{b_\xi^j, a_\xi\}) = 0 = c_n(\{b_\eta^j, a_\xi\})$, since $k \neq j$. Similarly, if $\zeta = \eta$, we have $c_n(\{b_\xi^j, a_\eta\}) = 0 = c_n(\{b_\eta^j, a_\eta\})$.

\medskip

In summary, $c_n(\{\alpha,\gamma\}) = c_n(\{\beta,\gamma\})$ for every $\gamma \in E_\alpha \cup E_\beta$, which proves that $c_n$ satisfies condition $(\star_n)$.
\end{proof}

%%%%%%%%%%%%%%%%%%%%%%%%%%%%%%%%%%%%%%%%%%%%%%%%%%%%%%%%%%%%%%
\section{Separation properties}\label{sec_SGT}
%%%%%%%%%%%%%%%%%%%%%%%%%%%%%%%%%%%%%%%%%%%%%%%%%%%%%%%%%%%%%%

In this section, we study the separation properties of the spaces $X_c$. We begin by characterizing when these spaces are Hausdorff and show that every coloring can be modified, while preserving certain combinatorial properties, to obtain Hausdorff spaces. We then prove that there exist Hausdorff spaces $X_c$ that are not normal and examine the separation properties that arise when $R(X_c) = \omega$.

\begin{lemma}\label{lema_SGT_Xc_T2} The following statements are equivalent for any coloring $c : [\omega_1]^2 \to 2$.

\begin{enumerate}
\item $X_c$ is $T_0$.
\item $X_c$ is $T_1$.
\item $X_c$ is $T_2$.
\item $X_c$ is $T_{3\frac{1}{2}}$.
\item For every $\alpha < \beta < \omega_1$, we have $c_\alpha \neq c_\beta$.
\item For every $\alpha < \beta < \omega_1$, we have $c(\{\alpha,\beta\}) = 1$ or $\Delta_c(\alpha,\beta) \neq \emptyset$.
\end{enumerate}

\end{lemma}

\begin{proof} The equivalence of (1)--(4) follows from the fact that $X_c$ is always zero-dimensional. Therefore, it suffices to show that (2) implies (5), that (5) implies (3), and that (5) and (6) are equivalent.

To prove that (2) implies (5), observe that if $\alpha < \beta < \omega_1$, then there exists $p \in \mathcal{P}_\alpha$ such that $\beta \notin U_p$; that is, there exists $\gamma \in \operatorname{dom}(p)$ such that $c_\alpha(\gamma) = p(\gamma) \neq c_\beta(\gamma)$. Consequently, $c_\alpha \neq c_\beta$.

On the other hand, to prove that (5) implies (3), suppose that (5) holds. Let $\alpha < \beta < \omega_1$ and $\gamma < \omega_1$ be such that $c_\alpha(\gamma) \neq c_\beta(\gamma)$. If $p,q : \{\gamma\} \to 2$ are defined by $p(\gamma) := c_\alpha(\gamma)$ and $q(\gamma) := c_\beta(\gamma) = 1 - c_\alpha(\gamma)$, then $p \in \mathcal{P}_\alpha$, $q \in \mathcal{P}_\beta$, and $U_p \cap U_q = \emptyset$.

Finally, the equivalence of (5) and (6) follows from the fact that, for every $\alpha < \beta < \omega_1$, the functions $c_\alpha$ and $c_\beta$ differ at some coordinate in $\{\alpha,\beta\}$ if and only if $c(\{\alpha,\beta\}) = 1$.
\end{proof}

Although not all spaces $X_c$ are Hausdorff, the following result shows that this property can be guaranteed by modifying the original coloring while preserving both $S_c$ and $S_{\overline{c}}$.

In Lemma~\ref{lema_SGT_coloracion_Hausdorff}, we adopt the convention that if $D \subseteq \omega_1$ and $t : D \to 2$ is a function, then $1-t : D \to 2$ denotes the function defined by $(1-t)(\gamma) := 1 - t(\gamma)$.

\begin{lemma}\label{lema_SGT_coloracion_Hausdorff} If $c : [\omega_1]^2 \to 2$ is a coloring and $D \in [\omega_1]^\omega$, then there exists a coloring $d : [\omega_1]^2 \to 2$ satisfying the following conditions:

\begin{enumerate}
\item $d(\{\alpha,\beta\}) = c(\{\alpha,\beta\})$ for every pair of distinct $\alpha,\beta \in D$ and every pair of distinct $\alpha,\beta \in \omega_1 \setminus D$.

\item $S_d = S_c$ and $S_{\overline{d}} = S_{\overline{c}}$.

\item The spaces $X_d$ and $X_{\overline{d}}$ are Hausdorff.

\end{enumerate}

\end{lemma}

\begin{proof} For each $\delta \in D$, let $s_\delta := \restr{c_\delta}{D}$ and $t_\delta := \restr{\overline{c}_\delta}{D}$. Clearly, the set
\[
\mathcal{E} :=
\{s_\delta : \delta \in D\}
\cup
\{1-t_\delta : \delta \in D\}
\]
is countable. Observe also that, for each $\delta \in D$, the functions $s_\delta$ and $1-t_\delta$ agree on $D\setminus\{\delta\}$, but $s_\delta(\delta) = 0$ and $(1-t_\delta)(\delta) = 1$.

\medskip

\noindent {\bf Claim.} There exists an injective function $h : \omega_1 \setminus D \to 2^D \setminus \mathcal{E}$ such that, for every pair of distinct $\delta,\varepsilon \in D$, there exists $\gamma \in \omega_1 \setminus D$ satisfying $h(\gamma)(\delta) \neq h(\gamma)(\varepsilon)$.

\medskip

First, since $|D| = \omega$, there exists an enumeration without repetitions of $[D]^2$, say $\{\{\delta_n,\varepsilon_n\} : n < \omega\}$, where $\delta_n \neq \varepsilon_n$ for every $n < \omega$. For each $n < \omega$, consider
\[
P_n = \big\{p \in 2^D : p(\delta_n) = 0 \ \wedge \ p(\varepsilon_n) = 1\big\}.
\]
Since $|P_n| = \mathfrak{c}$ and $|\mathcal{E}| \leq \omega$, we can recursively choose a family $\{p_n : n < \omega\}$ of pairwise distinct functions such that $p_n \in P_n \setminus \mathcal{E}$ for every $n < \omega$. In particular, $p_n(\delta_n) \neq p_n(\varepsilon_n)$ for every $n < \omega$.

Now, let $\{\gamma_n : n < \omega\}$ be a subset of $\omega_1 \setminus D$ enumerated without repetitions, and consider the function $h_0 : \{\gamma_n : n < \omega\} \to 2^D$ defined by $h_0(\gamma_n) := p_n$. On the other hand, since $2^D \setminus (\mathcal{E} \cup \{p_n : n < \omega\})$ has cardinality $\mathfrak{c}$, there exists an injective function
\[
h_1 : \left(\omega_1 \setminus D\right)
\setminus \{\gamma_n : n < \omega\}
\longrightarrow
2^D \setminus
\left(\mathcal{E} \cup \{p_n : n < \omega\}\right).
\]
Finally, the function $h := h_0 \cup h_1$ satisfies the desired properties.

\medskip

Let $d : [\omega_1]^2 \to 2$ be the coloring defined by
\[
d(\{\alpha,\beta\}) :=
\begin{cases}
h(\alpha)(\beta), & \text{if $\alpha \notin D$ and $\beta \in D$},\\
h(\beta)(\alpha), & \text{if $\alpha \in D$ and $\beta \notin D$},\\
c(\{\alpha,\beta\}), & \text{otherwise}.
\end{cases}
\]
Clearly, $d$ satisfies condition (1).

To prove condition (2), let $n \in \mathbb{N}$, let $J \in [\omega_1]^{\omega_1}$, and let $\{E_\alpha : \alpha \in J\}$ be an $n$-block sequence. Since this family is pairwise disjoint and $D$ is countable, there exists $I \in [J]^{\omega_1}$ such that $I \cap D = \emptyset$ and $E_\alpha \cap D = \emptyset$ for every $\alpha \in I$. Thus, for every pair of distinct $\alpha,\beta \in I$ and every $\xi \in \{\alpha,\beta\}$, we have $\restr{c_\xi}{E_\alpha \cup E_\beta} = \restr{d_\xi}{E_\alpha \cup E_\beta}$. This shows that $c$ satisfies condition $(\star_n^c)$ if and only if $d$ satisfies condition $(\star_n^d)$. The same argument, applied to the complementary colorings, shows that $\overline{c}$ satisfies $(\star_n^{\overline{c}})$ if and only if $\overline{d}$ satisfies $(\star_n^{\overline{d}})$. Therefore, $S_d = S_c$ and $S_{\overline{d}} = S_{\overline{c}}$.

Finally, to establish condition (3), we show that $\Delta_d(\alpha,\beta) \neq \emptyset$ for every pair of distinct $\alpha,\beta < \omega_1$. Let $\alpha,\beta < \omega_1$ be distinct. If $\alpha,\beta \notin D$, then $\restr{d_\alpha}{D} = h(\alpha)$ and $\restr{d_\beta}{D} = h(\beta)$. Since $h$ is injective, there exists $\gamma \in D$ such that $d_\alpha(\gamma) \neq d_\beta(\gamma)$. Therefore, $\gamma \in \Delta_d(\alpha,\beta)$.

Now, suppose that exactly one of the ordinals $\alpha,\beta$ belongs to $D$. Without loss of generality, we may assume that $\alpha \in D$ and $\beta \notin D$. Since the colorings $c$ and $d$ agree on $[D]^2$, we have $\restr{d_\alpha}{D} = s_\alpha$. On the other hand, $\restr{d_\beta}{D} = h(\beta)$. The choice of $h$ guarantees that $h(\beta) \neq s_\alpha$ and $h(\beta) \neq 1-t_\alpha$. Since $s_\alpha$ and $1-t_\alpha$ agree on $D\setminus\{\alpha\}$, there exists $\gamma \in D\setminus\{\alpha\}$ such that $h(\beta)(\gamma) \neq s_\alpha(\gamma)$. Consequently, $d_\alpha(\gamma) \neq d_\beta(\gamma)$, so $\gamma \in \Delta_d(\alpha,\beta)$.

Lastly, suppose that $\alpha,\beta \in D$. By the choice of $h$, there exists $\gamma \in \omega_1 \setminus D$ such that $h(\gamma)(\alpha) \neq h(\gamma)(\beta)$. Therefore, $d_\alpha(\gamma) \neq d_\beta(\gamma)$ and, consequently, $\gamma \in \Delta_d(\alpha,\beta)$.

Thus, $\Delta_d(\alpha,\beta) \neq \emptyset$ for every pair of distinct $\alpha,\beta < \omega_1$. Since $\Delta_d(\alpha,\beta) = \Delta_{\overline{d}}(\alpha,\beta)$, Lemma~\ref{lema_SGT_Xc_T2} implies that both $X_d$ and $X_{\overline{d}}$ are Hausdorff spaces.
\end{proof}

The spaces $X_c$ are zero-dimensional and therefore regular. However, this property does not guarantee that they are normal, even when they are Hausdorff. To show that this situation can occur within our class of spaces, we use the following auxiliary result.

\begin{lemma}\label{lema_SGT_representacion} If $Y$ is a zero-dimensional space such that $|Y| = \omega_1$ and $w(Y) \leq \omega_1$, then there exists a coloring $c : [\omega_1]^2 \to 2$ such that
\[
X_c \cong Y \oplus D(\omega_1),
\]
where $D(\omega_1)$ denotes the discrete space of cardinality $\omega_1$.

\end{lemma}

\begin{proof} Let $z \in Y$. By the assumptions on $Y$, there exists a family $\{C_\xi : \xi < \omega_1\}$ of clopen subsets of $Y$ such that $z \in \bigcap_{\xi < \omega_1} C_\xi$ and $\{C_\xi : \xi < \omega_1\} \cup \{Y \setminus C_\xi : \xi < \omega_1\}$ is a base for the topology of $Y$.

Let $A,B \in [\omega_1]^{\omega_1}$ be such that $A \cap B = \emptyset$ and $A \cup B = \omega_1$. Moreover, let $\{a_y : y \in Y\}$ and $\{b_\xi : \xi < \omega_1\}$ be enumerations without repetitions of $A$ and $B$, respectively. Also, let $c : [\omega_1]^2 \to 2$ be the coloring defined by
\[
c(\{\alpha,\beta\}) :=
\begin{cases}
0, & \text{if } \alpha,\beta \in A,\\
1, & \text{if } \alpha,\beta \in B,\\
1, & \text{if } \{\alpha,\beta\} = \{a_y,b_\xi\} \text{ with } y \in C_\xi,\\
0, & \text{if } \{\alpha,\beta\} = \{a_y,b_\xi\} \text{ with } y \notin C_\xi.
\end{cases}
\]

Observe that, since $z \in C_\xi$ for every $\xi < \omega_1$, we have $U_{a_z}^0 = A$ and $U_{a_z}^1 = B$. Therefore, $A$ and $B$ are clopen subsets of $X_c$. Moreover, for every $\xi < \omega_1$, we have $U_{b_\xi}^0 \cap B = \{b_\xi\}$. Consequently, $B$ is a discrete subspace of $X_c$ of cardinality $\omega_1$.

Now, let $f : Y \to A$ be the bijection defined by $f(y) := a_y$. By the definition of $c$, for every $\xi < \omega_1$ and $y \in Y$, we have
\[
U_{b_\xi}^1 \cap A = f[C_\xi], \hspace{0.5em} U_{b_\xi}^0 \cap A = f[Y \setminus C_\xi], \hspace{0.5em} U_{a_y}^0 \cap A = A \hspace{0.5em} \text{and} \hspace{0.5em} U_{a_y}^1 \cap A = \emptyset.
\]
Since the family $\{C_\xi : \xi < \omega_1\} \cup \{Y \setminus C_\xi : \xi < \omega_1\}$ is a base for the topology of $Y$, the preceding equalities imply that $f$ is a homeomorphism.

Finally, since $A$ and $B$ are disjoint open sets whose union is $X_c$, we conclude that $X_c \cong A \oplus B \cong Y \oplus D(\omega_1)$.
\end{proof}

\begin{proposition}\label{prop_SGT_no_normal} Let $Y$ be a zero-dimensional Hausdorff space that is not normal. If $|Y| = \omega_1$ and $w(Y) \leq \omega_1$, then there exists a coloring $c_Y : [\omega_1]^2 \to 2$ such that $X_{c_Y}$ is Hausdorff, zero-dimensional, and not normal.

\end{proposition}

\begin{proof} By Lemma~\ref{lema_SGT_representacion}, there exists a coloring $c_Y : [\omega_1]^2 \to 2$ such that $X_{c_Y} \cong Y \oplus D(\omega_1)$. Since $Y$ and $D(\omega_1)$ are Hausdorff spaces, so is $X_{c_Y}$. Moreover, $X_{c_Y}$ is zero-dimensional, as is every space associated with a coloring. Finally, since $Y$ is homeomorphic to a closed subspace of $X_{c_Y}$ and every closed subspace of a normal space is normal, we conclude that $X_{c_Y}$ is not normal.
\end{proof}

It is worth noting that the hypotheses of Proposition~\ref{prop_SGT_no_normal} can be satisfied in $\mathsf{ZFC}$. For example, the space
\[
Y := \big([0,\omega_1] \times [0,\omega]\big) \setminus \{(\omega_1,\omega)\},
\]
known as the \emph{deleted Tychonoff plank}, satisfies all of them.

On the other hand, if $X_c$ is $T_2$ and $R(X_c) = \omega$, then all its subspaces are $T_3$ and Lindelöf; consequently, they are also $T_6$, paracompact, and collectionwise normal.

Observe that the Hausdorff property can be guaranteed without losing the condition $R(X_c) = \omega$. Indeed, by Lemmas~\ref{lema_GT_Xc_HG} and \ref{lema_SGT_coloracion_Hausdorff}, for every coloring $c : [\omega_1]^2 \to 2$ such that $R(X_c) = \omega$, there exists a coloring $d : [\omega_1]^2 \to 2$ such that $R(X_d) = \omega$ and $X_d$ is Hausdorff.

Finally, note that if $\underline{1} : [\omega_1]^2 \to 2$ is the constant one coloring, then $X_{\underline{1}}$ is a discrete space of cardinality $\omega_1$ and, consequently, all its subspaces are $T_6$, paracompact, and collectionwise normal, but $R(X_{\underline{1}}) = \omega_1$.

%%%%%%%%%%%%%%%%%%%%%%%%%%%%%%%%%%%%%%%%%%%%%%%%%%%%%%%%%%%%%%
\section{Coloring principles}\label{sec_CP}
%%%%%%%%%%%%%%%%%%%%%%%%%%%%%%%%%%%%%%%%%%%%%%%%%%%%%%%%%%%%%%

The \emph{spread} of a topological space $X$ is the cardinal number
\[
s(X) := \sup\{|D| : \text{$D$ is a discrete subspace of $X$}\} + \omega.
\]
For our purposes, it suffices to recall that $s(X) \leq R(X)$.

We also introduce the following principles, which will be used throughout the remainder of this paper. Each asserts the existence of a coloring $c : [\omega_1]^2 \to 2$ satisfying the indicated condition:
\begin{align*}
(\mathsf{C}_1)&: \hspace{0.5em} c(\mathcal{F}[X_c]) = \omega = R(X_{\overline{c}}). \\
(\mathsf{C}_2)&: \hspace{0.5em} R(X_c) = \omega = R(X_{\overline{c}}). \\
(\mathsf{C}_3)&: \hspace{0.5em} R(X_c) = \omega < s(X_c^2). \\
(\mathsf{C}_4)&: \hspace{0.5em} R(X_c) = \omega < R(X_c^2).
\end{align*}

In view of the inequalities $s(X) \leq R(X) \leq c(\mathcal{F}[X])$, it is immediate that principle $(\mathsf{C}_1)$ implies $(\mathsf{C}_2)$ and that $(\mathsf{C}_3)$ implies $(\mathsf{C}_4)$. Next, we prove an auxiliary lemma that will allow us to establish, in Corollary~\ref{cor_CP_coloracion_d_e}, that $(\mathsf{C}_2)$ implies $(\mathsf{C}_3)$.

\begin{lemma}\label{lema_CP_coloracion_d_e} Let $d,e : [\omega_1]^2 \to 2$ be two colorings.

\begin{enumerate}
\item If $S_d = \mathbb{N} = S_e$ and $[\omega_1]^2 = d^{-1}\{1\} \cup e^{-1}\{1\}$, then there exists a coloring $c : [\omega_1]^2 \to 2$ such that $R(X_c) = \omega$ and $s(X_c^2) = \omega_1$.

\item Moreover, if the hypotheses of (1) hold and the spaces $X_d$ and $X_e$ are Hausdorff, then the coloring $c$ also satisfies that $X_c$ is Hausdorff.

\end{enumerate}

\end{lemma}

\begin{proof} Let $A,B \in [\omega_1]^{\omega_1}$ be such that $A \cup B = \omega_1$ and $A \cap B = \emptyset$. Moreover, let $\{a_\xi : \xi < \omega_1\}$ and $\{b_\xi : \xi < \omega_1\}$ be enumerations without repetitions of $A$ and $B$, respectively. Let $c : [\omega_1]^2 \to 2$ be the coloring defined by
\[
c(\{\alpha,\beta\}) :=
\begin{cases}
d(\{\xi,\eta\}), & \text{if $\{\alpha,\beta\} = \{a_\xi,a_\eta\}$ for some distinct $\xi,\eta < \omega_1$},\\
e(\{\xi,\eta\}), & \text{if $\{\alpha,\beta\} = \{b_\xi,b_\eta\}$ for some distinct $\xi,\eta < \omega_1$},\\
0, & \text{otherwise}.
\end{cases}
\]
The remainder of the proof is divided into a series of claims for clarity.

\medskip

\noindent {\bf Claim 1.} The subspaces $A$ and $B$ of $X_c$ are homeomorphic to $X_d$ and $X_e$, respectively.

\medskip

By symmetry, it suffices to prove the claim for one of the subspaces. Let $f : X_d \to A$ be the bijection defined by $f(\xi) := a_\xi$. Observe that, for every $\xi,\eta < \omega_1$, we have $c_{a_\xi}(a_\eta) = d_\xi(\eta)$ and $c_{a_\xi}(b_\eta) = 0$.

To prove that $f$ is continuous, let $p \in \mathrm{Fn}(\omega_1,2)$ and let $q : f^{-1}[\operatorname{dom}(p) \cap A] \to 2$ be defined by $q(\eta) := p(a_\eta)$. If there exists $\beta \in \operatorname{dom}(p) \cap B$ such that $p(\beta) = 1$, then $U_p(c) \cap A = \emptyset$, since $c_{a_\xi}(\beta) = 0$ for every $\xi < \omega_1$. Otherwise, the preceding identities imply that
\[
f^{-1}[U_p(c) \cap A] = U_q(d).
\]
Therefore, $f$ is continuous.

To prove that $f$ is open, let $q \in \mathrm{Fn}(\omega_1,2)$ and let $p : f[\operatorname{dom}(q)] \to 2$ be defined by $p(a_\eta) := q(\eta)$. Again, the preceding identities guarantee that
\[
f[U_q(d)] = U_p(c) \cap A.
\]
Thus, $f$ maps basic open sets of $X_d$ to open sets of the subspace $A$. Consequently, $f$ is open and hence a homeomorphism.

\medskip

Thus, Claim 1, together with the fact that $A \cup B = X_c$, implies that $R(X_c) = \omega$.

\medskip

\noindent {\bf Claim 2.} The set $D := \{(a_\xi,b_\xi) : \xi < \omega_1\}$ is a discrete subspace of $X_c^2$.

\medskip

For each $\xi < \omega_1$, let $U_\xi := U_{a_\xi}^{0}(c) \times U_{b_\xi}^{0}(c)$. Observe that $(a_\xi,b_\xi) \in U_\xi$ for every $\xi < \omega_1$. Moreover, if $\eta \in \omega_1 \setminus \{\xi\}$, then $(a_\eta,b_\eta) \in U_\xi$ if and only if $a_\eta \in U_{a_\xi}^{0}(c)$ and $b_\eta \in U_{b_\xi}^{0}(c)$; equivalently, if $c(\{a_\xi,a_\eta\}) = 0 = c(\{b_\xi,b_\eta\})$. The latter holds if and only if $d(\{\xi,\eta\}) = 0 = e(\{\xi,\eta\})$. However, this is impossible, since $\{\xi,\eta\} \in d^{-1}\{1\} \cup e^{-1}\{1\}$. Therefore, $U_\xi \cap D = \{(a_\xi,b_\xi)\}$ for every $\xi < \omega_1$.

\medskip

\noindent {\bf Claim 3.} If $X_d$ and $X_e$ are $T_2$, then $X_c$ is also $T_2$.

\medskip

By Lemma~\ref{lema_SGT_Xc_T2}, it suffices to show that the functions $c_\alpha$ are pairwise distinct. Let $\alpha,\beta < \omega_1$ be distinct.

If $\alpha,\beta \in A$, then there exist distinct $\xi,\eta < \omega_1$ such that $\alpha = a_\xi$ and $\beta = a_\eta$. Since $X_d$ is $T_2$, we have $d_\xi \neq d_\eta$. Consequently, $c_\alpha \neq c_\beta$. If $\alpha,\beta \in B$, a similar argument, using the fact that $X_e$ is $T_2$, yields the same conclusion.

On the other hand, suppose, without loss of generality, that $\alpha \in A$ and $\beta \in B$. Let $\xi,\eta < \omega_1$ be such that $\alpha = a_\xi$ and $\beta = b_\eta$. If $c_\alpha = c_\beta$, then both $d_\xi$ and $e_\eta$ would be constant zero functions, since $c_\beta$ vanishes on $A$ and $c_\alpha$ vanishes on $B$.

If $\xi \neq \eta$, this would imply that $d(\{\xi,\eta\}) = 0 = e(\{\xi,\eta\})$, contradicting the hypothesis $[\omega_1]^2 = d^{-1}\{1\} \cup e^{-1}\{1\}$. If $\xi = \eta$, it suffices to choose any $\zeta < \omega_1$ distinct from $\xi$ to obtain the same contradiction.

Thus, $c_\alpha \neq c_\beta$ in all cases, and by Lemma~\ref{lema_SGT_Xc_T2}, we conclude that $X_c$ is $T_2$.

\medskip

In summary, Claims 1, 2, and 3, together with the fact that $|X_c| = \omega_1$, imply that $R(X_c) = \omega$, $s(X_c^2) = \omega_1$, and $X_c$ is $T_2$ whenever $X_d$ and $X_e$ are.
\end{proof}

\begin{corollary}\label{cor_CP_coloracion_d_e} If principle $(\mathsf{C}_2)$ holds, then there exists a coloring $c : [\omega_1]^2 \to 2$ such that $X_c$ is Hausdorff, $R(X_c) = \omega$, and $s(X_c^2) = \omega_1$. In particular, principle $(\mathsf{C}_2)$ implies principle $(\mathsf{C}_3)$.

\end{corollary}

\begin{proof} Suppose that $(\mathsf{C}_2)$ holds, and let $d : [\omega_1]^2 \to 2$ be a coloring witnessing it. By Lemma~\ref{lema_SGT_coloracion_Hausdorff}, we may assume, without loss of generality, that $X_d$ and $X_{\overline{d}}$ are Hausdorff spaces. Thus, since $S_d = \mathbb{N} = S_{\overline{d}}$ and $[\omega_1]^2 = d^{-1}\{1\} \cup \overline{d}^{-1}\{1\}$, it suffices to apply Lemma~\ref{lema_CP_coloracion_d_e} to the colorings $d$ and $\overline{d}$.
\end{proof}

Regarding principle $(\mathsf{C}_1)$, our next goal is to show, in Corollary~\ref{cor_CP_SCCP_independencia}, that it is independent of $\mathsf{ZFC}$. We begin with the following auxiliary result:

\begin{lemma}\label{lema_CP_condicion_estrella_n} Let $n \in \mathbb{N}$ and let $c : [\omega_1]^2 \to 2$ be a coloring. If, for every $m < n+2$, every pairwise disjoint family $\{a_\alpha : \alpha < \omega_1\}$ of elements of $[\omega_1]^m$, and every function $M : m \times m \to 2$, there exist $\alpha < \beta < \omega_1$ such that
\[
c(\{a_\alpha(i),a_\beta(j)\}) = M(i,j)
\]
for every $i,j < m$, then $\{1,\ldots,n\} \subseteq S_c \cap S_{\overline{c}}$.
\end{lemma}

\begin{proof} By Lemma~\ref{lema_GT_estrella_baja}, to prove that $\{1,\ldots,n\} \subseteq S_c \cap S_{\overline{c}}$, it suffices to show that $c$ and $\overline{c}$ satisfy their respective conditions $(\star_n)$. Since the two proofs are similar, we carry out their initial steps simultaneously.

Let $J \in [\omega_1]^{\omega_1}$ and let $\{E_\alpha : \alpha \in J\}$ be an $n$-block sequence. For each $\alpha \in J$, let $T_\alpha := E_\alpha \cup \{\alpha\}$. By Lemma~\ref{lema_GT_adelgazamiento_2}, there exist $m \in \mathbb{N}$ with $m \leq n+1$, $r < m$, a family of colors $\{k_s : s \in m \setminus \{r\}\} \subseteq 2$, and a set $I \in [J]^{\omega_1}$ satisfying conditions (1), (2), and (3) of that lemma.

For each $i < 2$, let $M_i : m \times m \to 2$ be the function defined by
\[
M_i(s,t) :=
\begin{cases}
i, & \text{if } s = t = r,\\
k_s, & \text{if } t = r \text{ and } s \neq r,\\
k_t, & \text{if } s = r \text{ and } t \neq r,\\
0, & \text{otherwise}.
\end{cases}
\]
By hypothesis, for each $i < 2$, there exist $\alpha_i < \beta_i$ in $I$ such that
\[
c(\{t_{\alpha_i}^p,t_{\beta_i}^q\}) = M_i(p,q)
\]
for every $p,q < m$.

Observe that, for each $i < 2$, we have $c(\{\alpha_i,\beta_i\}) = c(\{t_{\alpha_i}^r,t_{\beta_i}^r\}) = M_i(r,r) = i$. On the other hand, if $s \in m \setminus \{r\}$, then $c(\{\alpha_i,t_{\alpha_i}^s\}) = k_s$ and $c(\{\beta_i,t_{\alpha_i}^s\}) = M_i(s,r) = k_s$. Therefore,
\[
c(\{\alpha_i,t_{\alpha_i}^s\}) = c(\{\beta_i,t_{\alpha_i}^s\}).
\]
Similarly,
\[
c(\{\alpha_i,t_{\beta_i}^s\}) = M_i(r,s) = k_s = c(\{\beta_i,t_{\beta_i}^s\}).
\]
Consequently,
\begin{equation}\label{eq_3}
c(\{\alpha_i,\gamma\}) = c(\{\beta_i,\gamma\})
\end{equation}
for every $i < 2$ and $\gamma \in (T_{\alpha_i} \cup T_{\beta_i}) \setminus \{\alpha_i,\beta_i\}$.

Now, if $\alpha_0 \in E_{\alpha_0}$, then $c(\{\alpha_0\}) = 0 = c(\{\beta_0,\alpha_0\})$. Similarly, if $\beta_0 \in E_{\beta_0}$, then $c(\{\alpha_0,\beta_0\}) = 0 = c(\{\beta_0\})$. Therefore, equality~\eqref{eq_3} implies that $c(\{\alpha_0,\gamma\}) = c(\{\beta_0,\gamma\})$ for every $\gamma \in E_{\alpha_0} \cup E_{\beta_0}$.

Finally, when $i = 1$, the equality $c(\{\alpha_1,\beta_1\}) = 1$ implies that $\overline{c}(\{\alpha_1,\beta_1\}) = 0$. Moreover, if $\gamma \in (E_{\alpha_1} \cup E_{\beta_1}) \setminus \{\alpha_1,\beta_1\}$, equality~\eqref{eq_3} implies that $\overline{c}(\{\alpha_1,\gamma\}) = \overline{c}(\{\beta_1,\gamma\})$. On the other hand, if $\alpha_1 \in E_{\alpha_1}$, then $\overline{c}(\{\alpha_1\}) = 0 = \overline{c}(\{\beta_1,\alpha_1\})$. Similarly, if $\beta_1 \in E_{\beta_1}$, then $\overline{c}(\{\alpha_1,\beta_1\}) = 0 = \overline{c}(\{\beta_1\})$. Consequently, $\overline{c}(\{\alpha_1,\gamma\}) = \overline{c}(\{\beta_1,\gamma\})$ for every $\gamma \in E_{\alpha_1} \cup E_{\beta_1}$.

In summary, $n \in S_c \cap S_{\overline{c}}$, which implies that $\{1,\ldots,n\} \subseteq S_c \cap S_{\overline{c}}$.
\end{proof}

To study the consistency of principle $(\mathsf{C}_1)$ and its relationship with the cardinal $\mathfrak{p}$, we begin with the following consequence of Theorem 3.1 of Kojman, Rinot, and Stepr\=ans \cite{krs2023}:

\begin{theorem}\label{thm_CP_KRS} The following statements are equivalent.

\begin{enumerate}[label=(\arabic*)]
\item $\operatorname{non}(\mathcal{M}) = \omega_1$.

\item There exists a coloring $c : [\omega_1]^2 \to \omega_1$ such that, for every $k,\ell < \omega$, every infinite pairwise disjoint family $\mathcal{A} \subseteq [\omega_1]^k$, and every uncountable pairwise disjoint family $\mathcal{B} \subseteq [\omega_1]^\ell$, there exists $a \in \mathcal{A}$ such that, for every function $M : k \times \ell \to \omega_1$, there exists $b \in \mathcal{B}$ satisfying
\[
c(\{a(i),b(j)\}) = M(i,j)
\]
for every $(i,j) \in k \times \ell$.

\end{enumerate}

\end{theorem}

In \cite[Definition~3(2)]{chenmertens2021}, Chen-Mertens, Kojman, and Stepr\=ans consider the principle $\operatorname{Pr}_0(\omega_1,\omega_1,\omega_1,\omega)$, which asserts the existence of a coloring $d : [\omega_1]^2 \to \omega_1$ such that, for every $n \in \mathbb{N}$, every uncountable pairwise disjoint family $\mathcal{F} \subseteq [\omega_1]^n$, and every function $M : n \times n \to \omega_1$, there exist $a,b \in \mathcal{F}$ with $\max a < \min b$ such that
\[
d(\{a(i),b(j)\}) = M(i,j)
\]
for every $i,j < n$.

For our purposes, it will suffice to work with this principle, whose formulation is simpler than the combinatorial condition in Theorem~\ref{thm_CP_KRS}. Our first goal is to establish a connection between the equality $\operatorname{non}(\mathcal{M}) = \omega_1$, the principle $\operatorname{Pr}_0(\omega_1,\omega_1,\omega_1,\omega)$, and principle $(\mathsf{C}_1)$.

\begin{proposition}\label{prop_CP_KRS_Pr0} If $\operatorname{non}(\mathcal{M}) = \omega_1$, then $\operatorname{Pr}_0(\omega_1,\omega_1,\omega_1,\omega)$ holds.

\end{proposition}

\begin{proof} Let $d : [\omega_1]^2 \to \omega_1$ be a coloring satisfying condition (2) of Theorem~\ref{thm_CP_KRS}. Moreover, let $n \in \mathbb{N}$, let $\mathcal{F} \subseteq [\omega_1]^n$ be an uncountable pairwise disjoint family, and let $M : n \times n \to \omega_1$ be a function. Choose $\mathcal{A} \in [\mathcal{F}]^\omega$ and fix $\delta < \omega_1$ such that $\bigcup\mathcal{A} \subseteq \delta$. Since $\mathcal{F}$ is pairwise disjoint, the family $\mathcal{B} := \{b \in \mathcal{F} : \min b > \delta\}$ is uncountable. By Theorem~\ref{thm_CP_KRS}, there exist $a \in \mathcal{A}$ and $b \in \mathcal{B}$ such that $d(\{a(i),b(j)\}) = M(i,j)$ for every $i,j < n$. Moreover, $\max a < \min b$, which establishes the desired principle.
\end{proof}

\begin{proposition}\label{prop_CP_non_implica_SCCP} The principle $\operatorname{Pr}_0(\omega_1,\omega_1,\omega_1,\omega)$ implies $(\mathsf{C}_1)$.

\end{proposition}

\begin{proof} Let $d : [\omega_1]^2 \to \omega_1$ be a coloring witnessing $\operatorname{Pr}_0(\omega_1,\omega_1,\omega_1,\omega)$. Consider the function $h : \omega_1 \to 2$ defined by
\[
h(\xi) :=
\begin{cases}
1, & \text{if } \xi = 1,\\
0, & \text{if } \xi \neq 1,
\end{cases}
\]
and let $c := h \circ d$.

Observe that, for every $m \in \mathbb{N}$, every $I \in [\omega_1]^{\omega_1}$, every pairwise disjoint family $\{a_\alpha : \alpha \in I\}$ of elements of $[\omega_1]^m$, and every function $M : m \times m \to 2$, there exist $\alpha,\beta \in I$ with $\alpha < \beta$ such that
\begin{equation}\label{eq_CP_matrices_SCCP}
c(\{a_\alpha(i),a_\beta(j)\}) = M(i,j)
\end{equation}
for every $i,j < m$.

Indeed, there exists $J \in [I]^{\omega_1}$ such that $\max a_\alpha < \min a_\beta$ for every $\alpha,\beta \in J$ with $\alpha < \beta$. Applying $\operatorname{Pr}_0(\omega_1,\omega_1,\omega_1,\omega)$ to the family $\{a_\alpha : \alpha \in J\}$ and to $M$, viewed as a function taking values in $\omega_1$, we obtain $\alpha,\beta \in J$ with $\alpha < \beta$ such that $d(\{a_\alpha(i),a_\beta(j)\}) = M(i,j)$ for every $i,j < m$. Since the restriction of $h$ to $2$ is the identity, equality~\eqref{eq_CP_matrices_SCCP} holds.

Thus, Lemma~\ref{lema_CP_condicion_estrella_n} guarantees that $S_c = \mathbb{N} = S_{\overline{c}}$ and, consequently, Lemma~\ref{lema_GT_Xc_HG} implies that $R(X_c) = \omega = R(X_{\overline{c}})$. Now, to prove that $c(\mathcal{F}[X_c]) = \omega$, we show that $R(X_c^\omega) = \omega$. By the relations~\eqref{eq_1}, it suffices to prove that $R(X_c^n) = \omega$ for every $n \in \mathbb{N}$.

Suppose, toward a contradiction, that there exists $n \in \mathbb{N}$ such that $R(X_c^n) > \omega$. Let $\{(x_\alpha,V_\alpha) : \alpha < \omega_1\}$ be an assignment in $X_c^n$ such that $\{V_\alpha : \alpha < \omega_1\}$ is a weak separation of the subspace $\{x_\alpha : \alpha < \omega_1\}$ of $X_c^n$.

Without loss of generality, we may assume that
\[
\textstyle
x_\alpha = (x_\alpha(i))_{i<n},
\quad \text{and} \quad
V_\alpha = \prod_{i<n}U_{p_\alpha^i}(c),
\]
where, for every $\alpha < \omega_1$ and $i < n$, there exists a set $F_\alpha^i \in [\omega_1]^{<\omega} \setminus \{\emptyset\}$ such that $p_\alpha^i = \restr{c_{x_\alpha(i)}}{F_\alpha^i}$.

For each $\alpha < \omega_1$, consider the finite set
\[
\textstyle
H_\alpha := \{x_\alpha(i) : i < n\} \cup \bigcup_{i<n}F_\alpha^i.
\]
By the $\Delta$-system lemma, there exist $J \in [\omega_1]^{\omega_1}$ and a finite set $R \subseteq \omega_1$ such that $H_\alpha \cap H_\beta = R$ for every pair of distinct $\alpha,\beta \in J$. For each $\alpha \in J$, let $A_\alpha := H_\alpha \setminus R$. Observe that the family $\{A_\alpha : \alpha \in J\}$ is pairwise disjoint.

Now, by an argument similar to the one used in the proof of Lemma~\ref{lema_GT_adelgazamiento_1}, we can find a set $I \in [J]^{\omega_1}$ and a natural number $m$ such that the following conditions hold:

\begin{enumerate}
\item $|A_\alpha| = m$ for every $\alpha \in I$. In what follows, $\{a_\alpha(j) : j < m\}$ denotes the strictly increasing enumeration of $A_\alpha$ for each $\alpha \in I$.

\item For each $i < n$, exactly one of the following alternatives holds:
\begin{enumerate}[label=(\alph*)]
\item There exists $r_i \in R$ such that $x_\alpha(i) = r_i$ for every $\alpha \in I$.
\item There exists $j_i < m$ such that $x_\alpha(i) = a_\alpha(j_i)$ for every $\alpha \in I$.
\end{enumerate}

\item For each $i < n$, there exist sets $D_i \subseteq R$ and $I_i \subseteq m$ such that
\[
\textstyle
F_\alpha^i = D_i \cup \{a_\alpha(j) : j \in I_i\}
\]
for every $\alpha \in I$.

\item There exists a function $M : m \times m \to 2$ such that $M(j,j) = 0$ for every $j < m$ and
\[
\textstyle
c(\{a_\alpha(j),a_\alpha(k)\}) = M(j,k)
\]
for every $\alpha \in I$ and every pair of distinct $j,k < m$.

\item There exists a function $N : R \times m \to 2$ such that
\[
\textstyle
c(\{r,a_\alpha(j)\}) = N(r,j)
\]
for every $\alpha \in I$, $r \in R$, and $j < m$.

\end{enumerate}

Observe that $M$ is symmetric, that is, $M(j,k) = M(k,j)$ for every $j,k < m$. Moreover, since the points $x_\alpha$ are pairwise distinct, there must be at least one coordinate $i < n$ satisfying alternative (2b). In particular, $m \geq 1$.

Since the family $\{A_\alpha : \alpha \in I\}$ is uncountable and pairwise disjoint, property~\eqref{eq_CP_matrices_SCCP} guarantees that there exist $\alpha,\beta \in I$ with $\alpha < \beta$ such that
\begin{equation}\label{eq_CP_matrices_cruzadas}
c(\{a_\alpha(j),a_\beta(k)\}) = M(j,k)
\end{equation}
for every $j,k < m$.

Our next goal is to show that $x_\alpha \in V_\beta$ and $x_\beta \in V_\alpha$. To this end, fix $i < n$ and distinguish the two alternatives in (2).

If alternative (2a) holds, there exists $r_i \in R$ such that $x_\alpha(i) = r_i = x_\beta(i)$. Since $x_\alpha \in V_\alpha$ and $x_\beta \in V_\beta$, it follows immediately that
\[
x_\alpha(i) \in U_{p_\beta^i}(c)
\quad\text{and}\quad
x_\beta(i) \in U_{p_\alpha^i}(c).
\]

On the other hand, if alternative (2b) holds, there exists $j_i < m$ such that $x_\alpha(i) = a_\alpha(j_i)$ and $x_\beta(i) = a_\beta(j_i)$. In this case, we show that
\[
\restr{c_{x_\alpha(i)}}{F_\beta^i} = p_\beta^i
\quad\text{and}\quad
\restr{c_{x_\beta(i)}}{F_\alpha^i} = p_\alpha^i.
\]

Let $\gamma \in F_\beta^i$. By (3), we have $F_\beta^i = D_i \cup \{a_\beta(j) : j \in I_i\}$. If $\gamma \in D_i$, then $\gamma \in R$ and, by (5),
\[
c(\{a_\alpha(j_i),\gamma\})
=
N(\gamma,j_i)
=
c(\{a_\beta(j_i),\gamma\}).
\]
On the other hand, if $\gamma = a_\beta(j)$ for some $j \in I_i$, then, by~\eqref{eq_CP_matrices_cruzadas} and (4),
\[
c(\{a_\alpha(j_i),a_\beta(j)\})
=
M(j_i,j)
=
c(\{a_\beta(j_i),a_\beta(j)\}).
\]
When $j = j_i$, the last equality follows from our convention that $c(\{\xi\}) = 0$ for every $\xi < \omega_1$, together with the condition $M(j_i,j_i) = 0$.

Thus, $\restr{c_{x_\alpha(i)}}{F_\beta^i} = \restr{c_{x_\beta(i)}}{F_\beta^i} = p_\beta^i$ and, therefore, $x_\alpha(i) \in U_{p_\beta^i}(c)$. Similarly, using the symmetry of $M$, we obtain $\restr{c_{x_\beta(i)}}{F_\alpha^i} = \restr{c_{x_\alpha(i)}}{F_\alpha^i} = p_\alpha^i$, which implies that $x_\beta(i) \in U_{p_\alpha^i}(c)$. Consequently,
\[
x_\alpha(i) \in U_{p_\beta^i}(c)
\quad\text{and}\quad
x_\beta(i) \in U_{p_\alpha^i}(c)
\]
for every $i < n$. Therefore, $x_\alpha \in V_\beta$ and $x_\beta \in V_\alpha$, contradicting the fact that the family $\{V_\gamma : \gamma < \omega_1\}$ is a weak separation.

Finally, the preceding argument implies that $R(X_c^n) = \omega$ for every $n \in \mathbb{N}$, which in turn guarantees that $R(X_c^\omega) = \omega$. Consequently, the inequalities~\eqref{eq_1} yield $c(\mathcal{F}[X_c]) = \omega$, establishing principle $(\mathsf{C}_1)$.
\end{proof}

Our next task is to show that principle $(\mathsf{C}_1)$ is incompatible with Martin's axiom for $\omega_1$:

\begin{proposition}\label{prop_CP_SCCP_implica_no_MA} If principle $(\mathsf{C}_1)$ holds, then $\mathsf{MA}(\omega_1)$ fails.

\end{proposition}

\begin{proof} Suppose, toward a contradiction, that both $(\mathsf{C}_1)$ and $\mathsf{MA}(\omega_1)$ hold. Let $c : [\omega_1]^2 \to 2$ be a coloring such that $c(\mathcal{F}[X_c]) = \omega = R(X_{\overline{c}})$. Moreover, for each $\alpha < \omega_1$, let $V_\alpha := [\{\alpha\},U_\alpha^0(c)]$. Observe that, since $\mathcal{F}[X_c]$ has countable cellularity and $\mathsf{MA}(\omega_1)$ holds, this space has precaliber $\omega_1$. Therefore, there exists a set $J \in [\omega_1]^{\omega_1}$ such that the family $\{V_\alpha : \alpha \in J\}$ is centered.

Now, since $V_\alpha \cap V_\beta \neq \emptyset$ for every pair of distinct $\alpha,\beta \in J$, it follows that $\beta \in U_\alpha^0(c)$ and, consequently, $c(\{\alpha,\beta\}) = 0$. Thus, $U_\alpha^0(\overline{c}) \cap J = \{\alpha\}$ for every $\alpha \in J$. Consequently, $J$ is an uncountable discrete subspace of $X_{\overline{c}}$, contradicting $R(X_{\overline{c}}) = \omega$.
\end{proof}

Since both $\operatorname{non}(\mathcal{M}) = \omega_1$ and $\mathsf{MA}(\omega_1)$ are relatively consistent with $\mathsf{ZFC}$, Propositions~\ref{prop_CP_non_implica_SCCP} and \ref{prop_CP_SCCP_implica_no_MA} yield the following corollary:

\begin{corollary}\label{cor_CP_SCCP_independencia} Principle $(\mathsf{C}_1)$ is independent of $\mathsf{ZFC}$.
\end{corollary}

By a theorem of Bell, the inequality $\mathfrak{p} > \omega_1$ is equivalent to Martin's axiom for $\omega_1$ restricted to $\sigma$-centered partial orders. Since every $\sigma$-centered partial order is $\sigma$-linked, we have
\[
\mathsf{MA}_{\omega_1}(\sigma\text{-linked})
\Longrightarrow
\mathsf{MA}_{\omega_1}(\sigma\text{-centered})
\Longleftrightarrow
\mathfrak{p} > \omega_1.
\]
On the other hand, Chen-Mertens, Kojman, and Stepr\=ans~\cite[Corollary~29]{chenmertens2021} proved that it is relatively consistent with $\mathsf{ZFC}$ that $\mathsf{MA}_{\omega_1}(\sigma\text{-linked})$ and $\operatorname{Pr}_0(\omega_1,\omega_1,\omega_1,\omega)$ hold simultaneously; in fact, their result establishes a stronger coloring property over partitions. Together with Proposition~\ref{prop_CP_non_implica_SCCP}, this yields the following consequence:

\begin{corollary}\label{cor_CP_SCCP_no_p} It is relatively consistent with $\mathsf{ZFC}$ that $(\mathsf{C}_1)$ and $\mathfrak{p} > \omega_1$ both hold. In particular, $(\mathsf{C}_1)$ does not imply that $\mathfrak{p} = \omega_1$.

\end{corollary}

Next, we examine whether $(\mathsf{C}_2)$ is compatible with $\mathsf{MA}(\omega_1)$. To this end, we use \emph{2-entangled} sets and a construction of Gruenhage~\cite[pp.~302--304]{gruenhage1989}.

Recall that a set $E = \{e_\xi : \xi < \omega_1\} \subseteq \mathbb{R}$ is \emph{2-entangled} if, for every $2$-block sequence $\{F_\alpha : \alpha < \omega_1\}$ and every function $\sigma : 2 \to \{<,>\}$, there exist $\alpha < \beta < \omega_1$ such that
\[
e_{F_\alpha(i)}\;\sigma(i)\;e_{F_\beta(i)}
\]
for every $i < 2$ (see~\cite[Definition~13 and Lemma~17]{guztod2024}).

In~\cite[p.~303]{gruenhage1989}, Gruenhage presents a construction based on $2$-entangled sets. We describe it below using notation convenient for our purposes.

Let $A = \{a_\alpha : \alpha < \omega_1\}$ be a subset of $E^2$, where $a_\alpha = (a_\alpha(0),a_\alpha(1))$, and suppose that $a_\alpha(i) = a_\beta(j)$ if and only if $\alpha = \beta$ and $i = j$, for every $\alpha,\beta < \omega_1$ and $i,j < 2$.

Now, for every $\alpha < \omega_1$ and $n < \omega$, let
\begin{align*}
D(\alpha,n) &:= \big\{\beta \in \omega_1 \setminus \{\alpha\} : \norm{a_\beta-a_\alpha} < 2^{-n}\big\}, \\
B_n^+(a_\alpha) &:= \{a_\alpha\} \cup
\big\{a_\beta : \beta \in D(\alpha,n) \ \wedge \
(a_\beta(0)-a_\alpha(0))(a_\beta(1)-a_\alpha(1)) > 0
\big\} \hspace{0.5em} \text{and} \\
B_n^-(a_\alpha) &:= \{a_\alpha\} \cup
\big\{a_\beta : \beta \in D(\alpha,n) \ \wedge \
(a_\beta(0)-a_\alpha(0))(a_\beta(1)-a_\alpha(1)) < 0
\big\},
\end{align*}
where $\norm{\cdot}$ denotes the usual Euclidean norm on $\mathbb{R}^{2}$. The families
\[
\mathcal{B}^+ := \{B_n^+(a_\alpha) : \alpha < \omega_1 \ \wedge \ n < \omega\}
\quad \text{and} \quad
\mathcal{B}^- := \{B_n^-(a_\alpha) : \alpha < \omega_1 \ \wedge \ n < \omega\}
\]
are bases for topologies $\tau^+$ and $\tau^-$ on $A$, respectively. Gruenhage shows that $(A,\tau^+)$ and $(A,\tau^-)$ are hereditarily Lindelöf and semimetric. Moreover, since every Lindelöf semimetric space has the HG property (see~\cite[p.~302]{gruenhage1989}), we have $R(A,\tau^+) = \omega = R(A,\tau^-)$.

With these results in place, we can relate $2$-entangled sets to principle $(\mathsf{C}_2)$:

\begin{proposition}\label{prop_CP_entangled_C2} The existence of a $2$-entangled set implies principle $(\mathsf{C}_2)$.

\end{proposition}

\begin{proof} Let $E \subseteq \mathbb{R}$ be a $2$-entangled set, let $(A,\tau^+)$ and $(A,\tau^-)$ be the topological spaces described above, and let $c : [\omega_1]^2 \to 2$ be the coloring defined by
\[
c(\{\alpha,\beta\}) :=
\begin{cases}
0, & \text{if } (a_\alpha(0)-a_\beta(0))(a_\alpha(1)-a_\beta(1)) > 0,\\
1, & \text{if } (a_\alpha(0)-a_\beta(0))(a_\alpha(1)-a_\beta(1)) < 0.
\end{cases}
\]
Our next goal is to show that the bijection $f : X_c \to (A,\tau^+)$ defined by $f(\alpha) := a_\alpha$ is open. To this end, it suffices to prove that $f$ maps subbasic open sets of $X_c$ to open sets of $(A,\tau^+)$.

Let $\alpha < \omega_1$, $i < 2$, and $a_\beta \in f[U_\alpha^i(c)]$. If $\beta = \alpha$, then $i = 0$ and $B_n^+(a_\alpha) \subseteq f[U_\alpha^0(c)]$ for every $n < \omega$. Thus, suppose that $\beta \neq \alpha$. If $n < \omega$ is such that
\[
2^{-n} <
\min\big\{
|a_\beta(0)-a_\alpha(0)|,\,
|a_\beta(1)-a_\alpha(1)|
\big\},
\]
then, for every $a_\gamma \in B_n^+(a_\beta)$ and $j < 2$, the numbers $a_\gamma(j)-a_\alpha(j)$ and $a_\beta(j)-a_\alpha(j)$ have the same sign. Consequently, $c(\{\alpha,\gamma\}) = c(\{\alpha,\beta\}) = i$, so $B_n^+(a_\beta) \subseteq f[U_\alpha^i(c)]$.

Thus, since $f$ is open and injective, it follows that $R(X_c) \leq R(A,\tau^+) = \omega$. Similarly, using the same function $f$ and the topology $\tau^-$, we obtain $R(X_{\overline{c}}) = \omega$. Consequently, principle $(\mathsf{C}_2)$ holds.
\end{proof}

On the other hand, Guzmán and Todor\v{c}evi\'c proved that, assuming the existence of a supercompact cardinal, it is consistent that $\mathsf{PID}$, $\mathsf{MA}(\omega_1)$, and the existence of a $2$-entangled set hold simultaneously (see~\cite[Theorem~97]{guztod2024}). This result, together with Proposition~\ref{prop_CP_entangled_C2}, yields the following consequence:

\begin{corollary}\label{cor_CP_C2_MA} It is consistent that $(\mathsf{C}_2)$, $\mathsf{PID}$, and $\mathsf{MA}(\omega_1)$ hold simultaneously. In particular, $(\mathsf{C}_2)$ implies neither $(\mathsf{C}_1)$ nor the failure of $\mathsf{MA}(\omega_1)$.

\end{corollary}

Although the preceding corollary separates principles $(\mathsf{C}_1)$ and $(\mathsf{C}_2)$, the question of whether $(\mathsf{C}_2)$ can be proved in $\mathsf{ZFC}$ remains open. It is worth noting, however, that we can obtain finite approximations to $(\mathsf{C}_2)$ without using additional axioms.

Indeed, to establish $(\mathsf{C}_2)$, we need a single coloring $c : [\omega_1]^2 \to 2$ such that $\{1,\ldots,n\} \subseteq S_c \cap S_{\overline{c}}$ for every $n \in \mathbb{N}$. Interchanging the universal and existential quantifiers yields a statement provable in $\mathsf{ZFC}$, as shown by the following proposition:

\begin{proposition}\label{prop_CP_aproximaciones_finitas} For each $n \in \mathbb{N}$, there exists a coloring $d_n : [\omega_1]^2 \to 2$ such that $\{1,\ldots,n\} \subseteq S_{d_n} \cap S_{\overline{d_n}}$.

\end{proposition}

\begin{proof} Let $n \in \mathbb{N}$. By Lemma~\ref{lema_GT_PW_n+1}, there exists a coloring $\widehat{c_n} : [\omega_1]^2 \to n+1$ such that, for every $m < n+2$, every pairwise disjoint family $\{a_\alpha : \alpha < \omega_1\}$ of elements of $[\omega_1]^m$, and every function $M : m \times m \to n+1$, there exist $\alpha < \beta < \omega_1$ such that
\[
\widehat{c_n}(\{a_\alpha(i),a_\beta(j)\}) = M(i,j)
\]
for every $i,j < m$.

If $h : n+1 \to 2$ is the function defined by
\[
h(i) :=
\begin{cases}
i, & \text{if } i < 2,\\
0, & \text{if } 2 \leq i < n+1,
\end{cases}
\]
and $d_n := h \circ \widehat{c_n}$, then $d_n$ satisfies all the hypotheses of Lemma~\ref{lema_CP_condicion_estrella_n}, so $\{1,\ldots,n\} \subseteq S_{d_n} \cap S_{\overline{d_n}}$.
\end{proof}

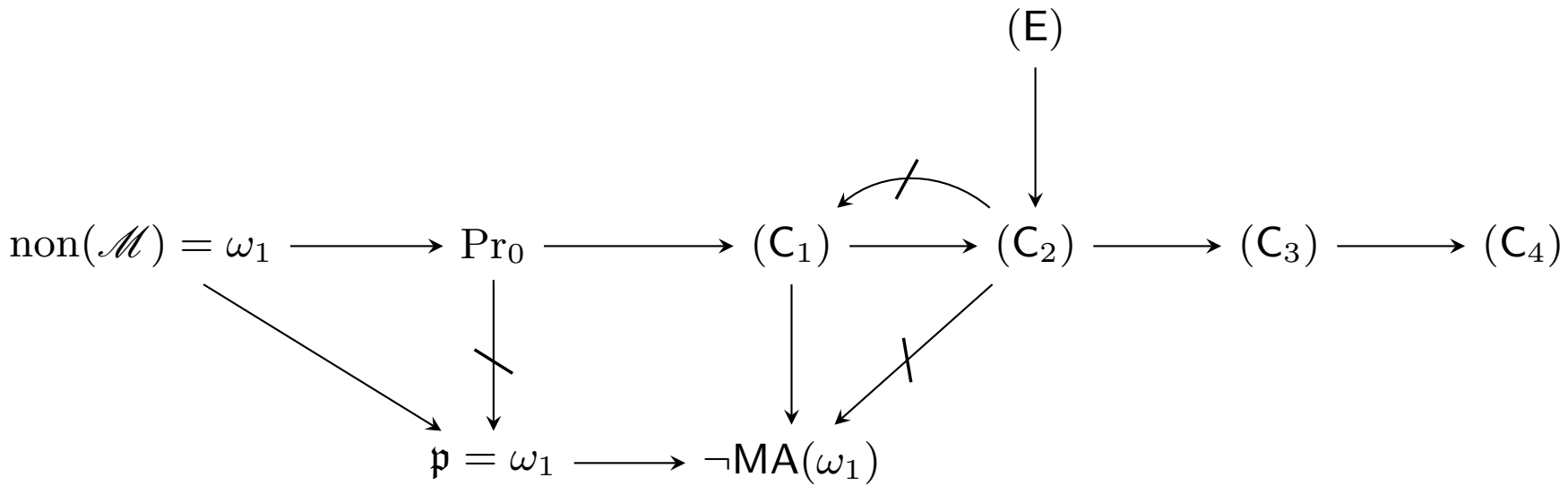
\begin{figure}[H]
\centering

\begin{tikzpicture}[
    >=stealth,
    every node/.style={font=\small},
    crossed/.style={
        postaction={decorate},
        decoration={
            markings,
            mark=at position 0.54 with {
                \draw[-,line width=0.65pt] (-2.5pt,-4pt) -- (2.5pt,4pt);
            }
        }
    }
]
    % First row
    \node (A) at (0,0) {$\operatorname{non}(\mathcal{M}) = \omega_1$};
    \node (B) at (2.6,0) {$\operatorname{Pr}_0$};
    \node (C) at (4.8,0) {$(\mathsf{C}_1)$};
    \node (D) at (6.6,0) {$(\mathsf{C}_2)$};
    \node (E) at (8.4,0) {$(\mathsf{C}_3)$};
    \node (F) at (10.2,0) {$(\mathsf{C}_4)$};

    % Second row and upper node
    \node (G) at (2.6,-1.6) {$\mathfrak{p} = \omega_1$};
    \node (H) at (4.8,-1.6) {$\neg\mathsf{MA}(\omega_1)$};
    \node (I) at (6.6,1.6) {$(\mathsf{E})$};

    % Implications
    \draw[->] (A) -- (B);
    \draw[->] (B) -- (C);
    \draw[->] (C) -- (D);
    \draw[->] (D) -- (E);
    \draw[->] (E) -- (F);
    \draw[->] (A) -- (G);
    \draw[->] (C) -- (H);
    \draw[->] (G) -- (H);
    \draw[->] (I) -- (D);

    % Non-implications
    \draw[->,crossed] (B) -- (G);
    \draw[->,crossed] (D)
        to[out=140,in=40] (C);
    \draw[->,crossed] (D) -- (H);
\end{tikzpicture}

\caption{Relations among combinatorial principles}
\label{fig_CP_implicaciones}
\end{figure}

Finally, the results of this section yield the following diagram of implications among the coloring principles and the set-theoretic statements considered here. In Figure~\ref{fig_CP_implicaciones}, we abbreviate $\operatorname{Pr}_0(\omega_1,\omega_1,\omega_1,\omega)$ by $\operatorname{Pr}_0$, and we write $(\mathsf{E})$ for the statement ``there exists a $2$-entangled set of real numbers.''

%%%%%%%%%%%%%%%%%%%%%%%%%%%%%%%%%%%%%%%%%%%%%%%%%%%%%%%%%%%%%%
\section{Pixley--Roy hyperspaces}\label{sec_PR}
%%%%%%%%%%%%%%%%%%%%%%%%%%%%%%%%%%%%%%%%%%%%%%%%%%%%%%%%%%%%%%

In this section, we study the relation between the weak separation number and the cellularity of Pixley--Roy hyperspaces. In particular, we prove that $c(\mathcal{F}[X_c]) = R(X_c^\omega)$ for any coloring $c : [\omega_1]^2 \to 2$ and that, if principle $(\mathsf{C}_2)$ holds, there exists a coloring $c : [\omega_1]^2 \to 2$ such that $X_c$ is $T_2$ and $R(X_c) < c(\mathcal{F}[X_c])$.

We begin with two results that hold for topological spaces in general.

\begin{lemma}\label{lema_PR_abiertos} If $X$ is a topological space, $n \in \mathbb{N}$, and $\{U_i : i < n\} \subseteq \tau_X$ is a pairwise disjoint family, then $R(\prod_{i<n} U_i) \leq c(\mathcal{F}[X])$.

\end{lemma}

\begin{proof} Let $A$ be a weakly separated subspace of $\prod_{i<n} U_i$, let $\kappa := |A|$, and let $\{x_\alpha : \alpha < \kappa\}$ be an enumeration without repetitions of $A$. For every $\alpha < \kappa$ and $i < n$, there exists an open set $V_\alpha(i)\subseteq U_i$ such that $x_\alpha(i) \in V_\alpha(i)$ and such that, defining $V_\alpha := \prod_{i<n} V_\alpha(i)$ for each $\alpha < \kappa$, the family $\{V_\alpha : \alpha < \kappa\}$ is a weak separation of $\{x_\alpha : \alpha < \kappa\}$.

Now, for each $\alpha < \kappa$, let $F_\alpha := \{x_\alpha(i) : i < n\}$ and $W_\alpha := \bigcup_{i<n} V_\alpha(i)$. Observe that $F_\alpha \subseteq W_\alpha$ for every $\alpha < \kappa$, so $[F_\alpha,W_\alpha]$ is a nonempty open subset of $\mathcal{F}[X]$. Moreover, if $\alpha < \beta < \kappa$ and we suppose that $[F_\alpha,W_\alpha] \cap [F_\beta,W_\beta] \neq \emptyset$, then $F_\alpha \subseteq W_\beta$ and $F_\beta \subseteq W_\alpha$. Since the family $\{U_i : i < n\}$ is pairwise disjoint, these inclusions imply that $x_\alpha(i) \in V_\beta(i)$ and $x_\beta(i) \in V_\alpha(i)$ for every $i < n$. Therefore, $x_\alpha \in V_\beta$ and $x_\beta \in V_\alpha$, contradicting the fact that $\{V_\gamma : \gamma < \kappa\}$ is a weak separation of $\{x_\gamma : \gamma < \kappa\}$. Consequently, the family $\{[F_\alpha,W_\alpha] : \alpha < \kappa\}$ is cellular and hence $\kappa \leq c(\mathcal{F}[X])$.
\end{proof}

\begin{proposition}\label{prop_PR_continua_inyectiva} Let $X$ be a topological space and $Y$ a second-countable Hausdorff space. If there exists a continuous injective function $f : X \to Y$, then $c(\mathcal{F}[X]) = c^*(\mathcal{F}[X]) = R(X^\omega)$.

\end{proposition}

\begin{proof} Let $\kappa := c(\mathcal{F}[X])$ and $\lambda := \kappa^+$. By the relations~\eqref{eq_1}, it suffices to prove that $R(X^n) \leq \kappa$ for every $n \in \mathbb{N}$. Suppose, toward a contradiction, that there exists $n \in \mathbb{N}$ such that $R(X^n) > \kappa$. Let $\{x_\alpha : \alpha < \lambda\}$ be an enumeration without repetitions of a weakly separated subspace of $X^n$. For each $\alpha < \lambda$, fix a family $\{V_\alpha(i) : i < n\} \subseteq \tau_X$ such that $x_\alpha(i) \in V_\alpha(i)$ for every $i < n$ and such that, defining $V_\alpha := \prod_{i<n} V_\alpha(i)$, the family $\{V_\alpha : \alpha < \lambda\}$ is a weak separation of $\{x_\alpha : \alpha < \lambda\}$.

For each $\alpha < \lambda$, let $\sim_\alpha$ be the equivalence relation on $n$ defined by $i \sim_\alpha j$ if and only if $x_\alpha(i) = x_\alpha(j)$. Since there are only finitely many equivalence relations on $n$, there exist an equivalence relation $\sim$ on $n$ and a set $I \in [\lambda]^\lambda$ such that $\sim_\alpha = \sim$ for every $\alpha \in I$. Let $m \in \mathbb{N}$ with $m \leq n$, and let $\{C_j : j < m\}$ be the set of equivalence classes of $\sim$.

For every $\alpha \in I$ and $j < m$, there exists a point $y_\alpha(j) \in X$ such that $x_\alpha(i) = y_\alpha(j)$ for every $i \in C_j$. Therefore, if for each $\alpha \in I$ we define the point $y_\alpha := (y_\alpha(0),\ldots,y_\alpha(m-1))$ of the product $X^m$, then the coordinates of each $y_\alpha$ are pairwise distinct.

Moreover, $\{y_\alpha : \alpha \in I\}$ is a weakly separated subspace of $X^m$. Indeed, for each $\alpha \in I$, consider
\[
\textstyle
W_\alpha :=
\prod_{j<m}
\left(
\bigcap_{i\in C_j} V_\alpha(i)
\right).
\]
The family $\{W_\alpha : \alpha \in I\}$ is a weak separation of $\{y_\alpha : \alpha \in I\}$, since, for every $\alpha,\beta \in I$, we have $y_\beta \in W_\alpha$ if and only if $x_\beta \in V_\alpha$.

Next, let $\mathcal{B}$ be a countable base for $Y$. Since $f$ is injective, for each $\alpha \in I$, the points in the family $\{f(y_\alpha(j)) : j < m\}$ are pairwise distinct. Thus, since $Y$ is Hausdorff, for each $\alpha \in I$, there exists a pairwise disjoint family $\{B_\alpha(j) : j < m\} \subseteq \mathcal{B}$ such that $f(y_\alpha(j)) \in B_\alpha(j)$ for every $j < m$.

Since $\mathcal{B}$ is countable and $m$ is finite, there exist a family $\{B_j : j < m\} \subseteq \mathcal{B}$ and a set $J \in [I]^\lambda$ such that $B_\alpha(j) = B_j$ for every $\alpha \in J$ and $j < m$. Now, for each $j < m$, let $U_j := f^{-1}[B_j]$. Clearly, $\{U_j : j < m\}$ is a pairwise disjoint family of open subsets of $X$. Moreover, by construction, we have $\{y_\alpha : \alpha \in J\} \subseteq \prod_{j<m} U_j$, where $\{y_\alpha : \alpha \in J\}$ is also a weakly separated subspace of $\prod_{j<m} U_j$.

Consequently, applying Lemma~\ref{lema_PR_abiertos}, we obtain the contradiction
\[
\textstyle
\kappa < \lambda \leq R(\prod_{j<m} U_j) \leq c(\mathcal{F}[X]) = \kappa.
\]
Therefore, $R(X^n) \leq c(\mathcal{F}[X])$ for every $n \in \mathbb{N}$.
\end{proof}

An equivalent formulation of the hypothesis of Proposition~\ref{prop_PR_continua_inyectiva} is that $X$ condenses onto a second-countable Hausdorff space. However, this is not equivalent to asserting that $X$ has countable $i$-weight, since the cardinal invariant $iw(X)$ is usually defined in the literature for Tychonoff spaces, whereas we do not require $X$ to satisfy this separation property.

It is also worth comparing Proposition~\ref{prop_PR_continua_inyectiva} with \cite[Corollary~4.14]{rios2023}, where it was proved that if $X$ is an infinite weakly separated space, then
\[
c^*(\mathcal{F}[X]) = R^*(X) = c(\mathcal{F}[X]) = R(X) = nw(X) = |X|.
\]
Thus, we obtain another class of spaces satisfying $c(\mathcal{F}[X]) = c^*(\mathcal{F}[X]) = R(X^\omega)$, without requiring weak separation.

Proposition~\ref{prop_PR_continua_inyectiva} will play a fundamental role in the proof of Theorem~\ref{thm_PR_equivalencia}, where we establish the equality of the cardinal invariants $c(\mathcal{F}[X_c])$, $c^*(\mathcal{F}[X_c])$, and $R(X_c^\omega)$ for every coloring $c : [\omega_1]^2 \to 2$. To this end, we introduce the following concept and study its relation to the spaces $X_c$:

\begin{definition}\label{def_PR_dep_coor} We say that a coloring $c : [\omega_1]^2 \to 2$ \emph{depends on fewer than $\omega_1$ coordinates} if there exists a set $D \in [\omega_1]^{\leq\omega}$ such that, for every $\alpha,\beta < \omega_1$, the equality $\restr{c_\alpha}{D} = \restr{c_\beta}{D}$ implies that $c_\alpha = c_\beta$; equivalently, if there exists a set $D \in [\omega_1]^{\leq\omega}$ such that the map $c_\alpha \mapsto \restr{c_\alpha}{D}$ is injective.

\end{definition}

\begin{lemma}\label{lema_PR_R_implica_dep} If $c : [\omega_1]^2 \to 2$ is a coloring and $R(X_c) = \omega$, then $c$ depends on fewer than $\omega_1$ coordinates.

\end{lemma}

\begin{proof} We proceed by contraposition. Suppose that $c$ does not depend on fewer than $\omega_1$ coordinates. This means that, for every set $D \in [\omega_1]^{\leq\omega}$, there exist $\alpha,\beta < \omega_1$ such that $\restr{c_\alpha}{D} = \restr{c_\beta}{D}$ but $c_\alpha \neq c_\beta$.

By transfinite recursion, we can construct families $\{a_\xi : \xi < \omega_1\}$, $\{b_\xi : \xi < \omega_1\}$, and $\{\gamma_\xi : \xi < \omega_1\}$ of ordinals below $\omega_1$ such that, for each $\xi < \omega_1$, we have $c_{a_\xi}(\gamma_\xi) \neq c_{b_\xi}(\gamma_\xi)$ and $c_{a_\xi}(\gamma_\eta) = c_{b_\xi}(\gamma_\eta)$ for every $\eta < \xi$.

Now, for each $\xi < \omega_1$, define
\[
\textstyle
U_\xi :=
U_{a_\xi}^{c_{a_\xi}(\gamma_\xi)}
\cap
U_{b_\xi}^{c_{b_\xi}(\gamma_\xi)}.
\]
Clearly, $U_\xi$ is an open subset of $X_c$ such that $\gamma_\xi \in U_\xi$. Moreover, the ordinals $\gamma_\xi$ are pairwise distinct, since if $\eta < \xi$ and $\gamma_\eta = \gamma_\xi$, we would have both $c_{a_\xi}(\gamma_\eta) = c_{b_\xi}(\gamma_\eta)$ and $c_{a_\xi}(\gamma_\xi) \neq c_{b_\xi}(\gamma_\xi)$.

On the other hand, if $\eta < \xi < \omega_1$, then $c_{a_\xi}(\gamma_\eta) = c_{b_\xi}(\gamma_\eta)$ and, consequently, the condition $c_{a_\xi}(\gamma_\xi) \neq c_{b_\xi}(\gamma_\xi)$ prevents $\gamma_\eta$ from belonging to $U_\xi$. Thus, the family $\{U_\xi : \xi < \omega_1\}$ is a weak separation of the subspace $\{\gamma_\xi : \xi < \omega_1\}$ of $X_c$, so $R(X_c) > \omega$.
\end{proof}

Although not directly related to the main goal of this section, the preceding lemma has an interesting topological consequence worth mentioning.

\begin{corollary}\label{cor_PR_no_pseudocompactos} If $X_c$ is $T_2$ and $R(X_c) = \omega$, then $X_c$ is not pseudocompact.

\end{corollary}

\begin{proof} Suppose, toward a contradiction, that $X_c$ is pseudocompact. As explained at the end of Section~\ref{sec_SGT}, $X_c$ is also paracompact and hence compact. By Lemma~\ref{lema_PR_R_implica_dep}, there exists a set $D \in [\omega_1]^{\leq\omega}$ such that the function $f : X_c \to 2^D$ defined by $f(\alpha) := \restr{c_\alpha}{D}$ is continuous. Moreover, by Lemma~\ref{lema_SGT_Xc_T2}, this function is injective. Since $X_c$ is compact and $2^D$ is Hausdorff, $f$ is a topological embedding. Thus, $X_c$ is second-countable and, being compact, has at most countably many clopen subsets. This contradicts the fact that the family $\{U_\alpha^1 : \alpha < \omega_1\}$ consists of pairwise distinct sets.
\end{proof}

The Hausdorff hypothesis in Corollary~\ref{cor_PR_no_pseudocompactos} is essential. Indeed, if $\underline{0} : [\omega_1]^2 \to 2$ is the constant zero coloring, then $X_{\underline{0}}$ is an indiscrete space and hence compact, but $R(X_{\underline{0}}) = \omega$.

On the other hand, there also exist compact Hausdorff spaces $X_c$. To see this, fix a partition $\mathcal{P}$ of $\omega_1 \setminus \{0\}$ into two-element sets and consider the coloring $c : [\omega_1]^2 \to 2$ defined by
\[
c(\{\alpha,\beta\}) :=
\begin{cases}
1, & \text{if } \{\alpha,\beta\} \in \mathcal{P},\\
0, & \text{if } \{\alpha,\beta\} \notin \mathcal{P}.
\end{cases}
\]
In this case, $X_c$ is homeomorphic to the Alexandroff compactification of the discrete space of cardinality $\omega_1$. In particular, $X_c$ is compact and Hausdorff, but $R(X_c) = \omega_1$.

We now return to our goal of establishing the equality of the cardinal invariants $c(\mathcal{F}[X_c])$, $c^*(\mathcal{F}[X_c])$, and $R(X_c^\omega)$. The following lemma shows how dependence on fewer than $\omega_1$ coordinates allows us to apply Proposition~\ref{prop_PR_continua_inyectiva}.

\begin{lemma}\label{lema_PR_dep_implica_c_c*_R} If a coloring $c : [\omega_1]^2 \to 2$ depends on fewer than $\omega_1$ coordinates, then $c(\mathcal{F}[X_c]) = c^*(\mathcal{F}[X_c]) = R(X_c^\omega)$.

\end{lemma}

\begin{proof} Let $D \in [\omega_1]^{\leq \omega}$ be a set witnessing that $c$ depends on fewer than $\omega_1$ coordinates. Observe that, if $Q_c := \{c_\alpha : \alpha < \omega_1\} \subseteq 2^{\omega_1}$, then the map $Q_c \to 2^D$ defined by $c_\alpha \mapsto \restr{c_\alpha}{D}$ is continuous and injective. Moreover, since $D$ is countable, the space $2^D$ is second-countable and Hausdorff. Therefore, Proposition~\ref{prop_PR_continua_inyectiva} guarantees that $c(\mathcal{F}[Q_c]) = c^*(\mathcal{F}[Q_c]) = R(Q_c^\omega)$.

For our purposes, since $c(\mathcal{F}[X_c]) \leq c^*(\mathcal{F}[X_c]) = R(X_c^\omega)$ (see~\eqref{eq_1}), it suffices to show that $R(X_c^\omega) \leq R(Q_c^\omega)$ and $c(\mathcal{F}[Q_c]) \leq c(\mathcal{F}[X_c])$.

To this end, consider the function $q : X_c \to Q_c$ defined by $q(\alpha) := c_\alpha$. Observe that, if for each $p \in \mathrm{Fn}(\omega_1,2)$ we set $B_p := \{x \in Q_c : p \subseteq x\}$, then $\{B_p : p \in \mathrm{Fn}(\omega_1,2)\}$ is a base for $Q_c$. Moreover, since $U_p = q^{-1}[B_p]$ for every $p \in \mathrm{Fn}(\omega_1,2)$, the function $q$ is continuous, open, and surjective. Furthermore, since $\{U_p : p \in \mathrm{Fn}(\omega_1,2)\}$ is a base for $X_c$, it follows that $\tau_{X_c} = \{q^{-1}[W] : W \in \tau_{Q_c}\}$ and, consequently, that $\tau_{X_c^\omega} = \{(q^\omega)^{-1}[W] : W \in \tau_{Q_c^\omega}\}$.

To prove that $R(X_c^\omega) \leq R(Q_c^\omega)$, let $\kappa$ be a cardinal and let $\{(x_\alpha,U_\alpha) : \alpha < \kappa\}$ be an assignment such that $\{U_\alpha : \alpha < \kappa\}$ is a weak separation of the subspace $\{x_\alpha : \alpha < \kappa\}$ of $X_c^\omega$. For each $\alpha < \kappa$, fix $W_\alpha \in \tau_{Q_c^\omega}$ such that $U_\alpha = (q^\omega)^{-1}[W_\alpha]$. Under these assumptions, the family $\{W_\alpha : \alpha < \kappa\}$ is a weak separation of the subspace $\{q^\omega(x_\alpha) : \alpha < \kappa\}$ of $Q_c^\omega$, since, for every $\alpha,\beta < \kappa$, we have $x_\alpha \in U_\beta$ if and only if $q^\omega(x_\alpha) \in W_\beta$. Thus, since $\{q^\omega(x_\alpha) : \alpha < \kappa\}$ has cardinality $\kappa$, it follows that $\kappa \leq R(Q_c^\omega)$.

Finally, to show that $c(\mathcal{F}[Q_c]) \leq c(\mathcal{F}[X_c])$, let $\{[F_\alpha,W_\alpha] : \alpha < \kappa\}$ be a cellular family of nonempty basic open sets in $\mathcal{F}[Q_c]$. For each $\alpha < \kappa$, choose $E_\alpha \in \mathcal{F}(X_c)$ with $q[E_\alpha] = F_\alpha$, and set $U_\alpha := q^{-1}[W_\alpha]$. Then $\{[E_\alpha,U_\alpha] : \alpha < \kappa\}$ is a cellular family in $\mathcal{F}[X_c]$, so $\kappa \leq c(\mathcal{F}[X_c])$.

Consequently, the equality $c(\mathcal{F}[Q_c]) = c^*(\mathcal{F}[Q_c]) = R(Q_c^\omega)$, together with the inequalities $R(X_c^\omega) \leq R(Q_c^\omega)$ and $c(\mathcal{F}[Q_c]) \leq c(\mathcal{F}[X_c])$, implies that $c(\mathcal{F}[X_c]) = c^*(\mathcal{F}[X_c]) = R(X_c^\omega)$.
\end{proof}

\begin{theorem}\label{thm_PR_equivalencia} If $c : [\omega_1]^2 \to 2$ is a coloring, then
\[
c(\mathcal{F}[X_c]) = c^*(\mathcal{F}[X_c]) = R(X_c^\omega).
\]

\end{theorem}

\begin{proof} First, observe that, since $|X_c| = \omega_1$, the cardinal $R(X_c)$ can only be $\omega$ or $\omega_1$. On the one hand, if $R(X_c) = \omega$, then Lemma~\ref{lema_PR_R_implica_dep} guarantees that $c$ depends on fewer than $\omega_1$ coordinates and, consequently, Lemma~\ref{lema_PR_dep_implica_c_c*_R} implies that $c(\mathcal{F}[X_c]) = c^*(\mathcal{F}[X_c]) = R(X_c^\omega)$. On the other hand, if $R(X_c) = \omega_1$, then, by~\eqref{eq_1}, we have
\begin{align*}
\omega_1 = R(X_c) \leq c(\mathcal{F}[X_c]) &\leq c^*(\mathcal{F}[X_c]) \\
&= R(X_c^\omega) = R^*(X_c) \leq nw(X_c) \leq |X_c| = \omega_1,
\end{align*}
which again yields $c(\mathcal{F}[X_c]) = c^*(\mathcal{F}[X_c]) = R(X_c^\omega)$.
\end{proof}

In particular, Theorem~\ref{thm_PR_equivalencia} shows that principle $(\mathsf{C}_1)$ is equivalent to the existence of a coloring $c : [\omega_1]^2 \to 2$ such that $R(X_c^\omega) = \omega = R(X_{\overline{c}})$.

In \cite[Question~4.16]{rios2023}, the author raised the question of whether there exists a space $X$ such that $R(X) < c(\mathcal{F}[X])$, or a space $Y$ such that $c(\mathcal{F}[Y]) < c^*(\mathcal{F}[Y])$. In view of Theorem~\ref{thm_PR_equivalencia}, the spaces $X_c$, which are the main objects of study in this paper, cannot provide an affirmative answer to the second part of this question. However, the results presented in this paper yield a consistent affirmative answer to the first part:

\begin{theorem}\label{thm_PR_HG_cel_no_num} If principle $(\mathsf{C}_2)$ holds, then there exists a coloring $c : [\omega_1]^2 \to 2$ such that $X_c$ is $T_2$, $R(X_c) = \omega$, and $c(\mathcal{F}[X_c]) = \omega_1$.

\end{theorem}

\begin{proof} By Corollary~\ref{cor_CP_coloracion_d_e}, there exists a coloring $c : [\omega_1]^2 \to 2$ such that $X_c$ is $T_2$, $R(X_c) = \omega$, and $s(X_c^2) = \omega_1$. Now, by Theorem~\ref{thm_PR_equivalencia}, we have
\[
\omega_1 = s(X_c^2) \leq R(X_c^2) \leq R(X_c^\omega) = c(\mathcal{F}[X_c]) \leq \omega_1.
\]
Therefore, $c(\mathcal{F}[X_c]) = \omega_1$.
\end{proof}

It is worth noting that the space $X_c$ in Theorem~\ref{thm_PR_HG_cel_no_num} also satisfies $R(X_c) = \omega < R(X_c^2)$, since Corollary~\ref{cor_CP_coloracion_d_e} guarantees that $s(X_c^2) = \omega_1$.

On the other hand, since $\mathsf{CH}$ implies $\operatorname{non}(\mathcal{M}) = \omega_1$, and the latter implies principle $(\mathsf{C}_2)$, Theorem~\ref{thm_PR_HG_cel_no_num} has two interesting consequences. First, as already mentioned, it provides a consistent affirmative answer to the first part of the question posed in \cite[Question~4.16]{rios2023}. Second, it strengthens a result of Hajnal and Juhász \cite[Theorem~1]{hajju1979}, who proved that, under $\mathsf{CH}$, there exists a coloring $c : [\omega_1]^2 \to 2$ such that $X_c$ is $T_2$, $R(X_c) = \omega$, and $nw(X_c) = \omega_1$.

%%%%%%%%%%%%%%%%%%%%%%%%%%%%%%%%%%%%%%%%%%%%%%%%%%%%%%%%%%%%%%
\section{Questions}\label{sec_Q}
%%%%%%%%%%%%%%%%%%%%%%%%%%%%%%%%%%%%%%%%%%%%%%%%%%%%%%%%%%%%%%

We conclude the paper with some open questions related to the problems studied here, including those posed in \cite[Question~4.16]{rios2023} and \cite[Question~3, p.~22]{fengar2026}.

\begin{question}\label{Q_1} Is there an additional axiom, consistent with $\mathsf{ZFC}$, under which the property $R(\boldsymbol{\cdot}) = \omega$ is finitely productive?

\end{question}

\begin{question}\label{Q_2} Does there exist in $\mathsf{ZFC}$ a topological space $X$ such that $R(X) < R(X^2)$?

\end{question}

\begin{question}\label{Q_3} Does there exist in $\mathsf{ZFC}$ a topological space $X$ such that $R(X) < c(\mathcal{F}[X])$?

\end{question}

\begin{question}\label{Q_4} Under what additional axiom, consistent with $\mathsf{ZFC}$, does there exist a topological space $X$ such that $c(\mathcal{F}[X]) = \omega < c^*(\mathcal{F}[X])$?

\end{question}

\begin{question}\label{Q_5} Is it true in $\mathsf{ZFC}$ that if $c(\mathcal{F}[X]) = \omega$, then $R(X^\omega) = \omega$?

\end{question}

\begin{question}\label{Q_6} Which of the principles $(\mathsf{C}_2)$, $(\mathsf{C}_3)$, and $(\mathsf{C}_4)$ are provable in $\mathsf{ZFC}$?
\end{question}

\begin{question}\label{Q_7} Are principles $(\mathsf{C}_2)$ and $(\mathsf{C}_3)$ equivalent? What about principles $(\mathsf{C}_3)$ and $(\mathsf{C}_4)$?

\end{question}

\Needspace{14\baselineskip}
\section*{AI Declaration}

During the development and preparation of this paper, the author used ChatGPT (OpenAI, GPT-6) to assist in exploring mathematical ideas, formulating results, developing and examining arguments, identifying bibliographic references, and improving the organization and exposition of the manuscript. The process involved extensive discussions and revisions, during which suggestions provided by ChatGPT were critically evaluated, corrected, or discarded when necessary, particularly in cases involving erroneous or incomplete arguments. All mathematical statements and proofs were carefully reviewed by the author, who assumes full responsibility for the content of the paper, including any remaining errors.


\begin{thebibliography}{99}

  \bibitem{chenmertens2021} W. Chen-Mertens, M. Kojman, J. Stepr\=ans, \emph{Strong colorings over partitions}, Bull. Symb. Log., {\bf 27}(1) (2021), 67--90, \url{https://doi.org/10.1017/bsl.2021.5}.

  \bibitem{fengar2026} Z. Feng, P. Gartside, \emph{The Base and Network Structures I: Chain Conditions and Tukey Types}, to appear in Acta Math. Hungar.

  \bibitem{gruenhage1989} G. Gruenhage, \emph{Cosmicity of cometrizable spaces}, Trans. Amer. Math. Soc., {\bf 313}(1) (1989), 301--315, \url{https://doi.org/10.1090/S0002-9947-1989-0992600-5}.

  \bibitem{guztod2024} O. Guzm\'an, S. Todor\v{c}evi\'c, \emph{The P-ideal dichotomy, Martin's axiom and entangled sets}, Israel J. Math., {\bf 263} (2024), 909--963, \url{https://doi.org/10.1007/s11856-024-2651-8}.

  \bibitem{hajju1979} A. Hajnal, I. Juhász, \emph{Weakly separated subspaces and networks}, Logic Colloquium '78 (Mons, 1978), North-Holland, Amsterdam, pp. 235--245, 1979.

  \bibitem{harkun2020} J. E. Hart, K. Kunen, \emph{Properties and super properties}, Topol. Proc., {\bf 55} (2020), 175--185.

  \bibitem{jss1996} I. Juhász, L. Soukup, Z. Szentmiklóssy, \emph{Forcing countable networks for spaces satisfying $R(X^\omega)=\omega$}, Comment. Math. Univ. Carolin., {\bf 37}(1) (1996), 159--170.

  \bibitem{krs2023} M. Kojman, A. Rinot, J. Stepr\=ans, \emph{Ramsey theory over partitions III: Strongly Luzin sets and partition relations}, Proc. Amer. Math. Soc., {\bf 151} (2023), 369--384, \url{https://doi.org/10.1090/proc/16106}.

  \bibitem{penwu2018} Y. Peng, L. Wu, \emph{A Lindel\"of group with non-Lindel\"of square}, Adv. Math., {\bf 325} (2018), 215--242, \url{https://doi.org/10.1016/j.aim.2017.11.021}.

  \bibitem{rios2023} A. Ríos Herrejón, \emph{Weakly separated spaces and Pixley--Roy hyperspaces}, Bull. Iran. Math. Soc., {\bf 49} (2023), 58, \url{https://doi.org/10.1007/s41980-023-00800-9}.

  \bibitem{tkachenko1978} M. G. Tkachenko, \emph{Chains and cardinals}, Dokl. Akad. Nauk SSSR, {\bf 239}(3) (1978), 546--549; English translation: Soviet Math. Dokl., {\bf 19}(2) (1978), 382--385.

\end{thebibliography}
\end{document}